\documentclass{amsart}
\usepackage{verbatim,amsmath,amssymb,amscd,color,graphics}

\usepackage[all]{xy} 
\CompileMatrices

\newcommand{\aal}{\mathbf{\stackrel{\leftarrow}{a}}}
\newcommand{\aar}{\mathbf{\stackrel{\rightarrow}{a}}}
\newcommand{\ab}{\ol{1}}
\newcommand{\bb}{\ol{2}}
\newcommand{\bij}{\Phi}

\newcommand{\eh}{\hat{e}}

\newcommand{\fh}{\hat{f}}

\newcommand{\geh}{\mathfrak{g}}
\newcommand{\id}{\mathrm{id}}
\newcommand{\Ih}{\hat{I}}

\newcommand{\is}{\mathrm{inner}}

\newcommand{\la}{\lambda}
\newcommand{\La}{\Lambda}

\newcommand{\ol}{\overline}

\newcommand{\os}{\mathrm{outer}}
\newcommand{\ot}{\otimes}
\newcommand{\pr}{\mathrm{pr}}

\newcommand{\sigD}{\mathfrak{S}}
\newcommand\tab[1]{\begin{array}{|c|}\hline {\lower1pt\hbox{$#1$}}\\ \hline \end{array}}
\newcommand{\Vh}{\hat{V}}
\newcommand{\ve}{\varepsilon}
\newcommand{\vp}{\varphi}
\newcommand{\wt}{\mathrm{wt}\,}
\newcommand{\Z}{\mathbb{Z}}

\newcommand{\hdom}{{\begin{picture}(8,4)\multiput(0,0)(4,0){3}{\line(0,1){4}}
\multiput(0,0)(0,4){2}{\line(1,0){8}}\end{picture}}}
\newcommand{\vdom}{{\begin{picture}(4,8)\multiput(0,0)(0,4){3}{\line(1,0){4}}
\multiput(0,0)(4,0){2}{\line(0,1){8}}\end{picture}}}
\newcommand{\cell}{{\begin{picture}(4,4)\multiput(0,0)(4,0){2}{\line(0,1){4}}
\multiput(0,0)(0,4){2}{\line(1,0){4}}\end{picture}}}

\numberwithin{equation}{section}

\newtheorem{theorem}{Theorem}
\newtheorem{proposition}[theorem]{Proposition}
\newtheorem{lemma}[theorem]{Lemma}
\newtheorem{corollary}[theorem]{Corollary}

\theoremstyle{definition}

\newtheorem{definition}{Definition}
\newtheorem{remark}{Remark}

\numberwithin{theorem}{section}
\numberwithin{definition}{section}
\numberwithin{remark}{section}

\begin{document}

\title[KR crystals]
{Kirillov--Reshetikhin crystals for nonexceptional types}

\author[G.~Fourier]{Ghislain Fourier}
\address{Mathematisches Institut der Universit\"at zu K\"oln,
Weyertal 86-90, 50931 K\"oln, Germany}
\email{gfourier@mi.uni-koeln.de}

\author[M.~Okado]{Masato Okado}
\address{Department of Mathematical Science,
Graduate School of Engineering Science, Osaka University,
Toyonaka, Osaka 560-8531, Japan}
\email{okado@sigmath.es.osaka-u.ac.jp}

\author[A.~Schilling]{Anne Schilling}
\address{Department of Mathematics, University of California, One Shields
Avenue, Davis, CA 95616-8633, U.S.A.}
\email{anne@math.ucdavis.edu}
\urladdr{http://www.math.ucdavis.edu/\~{}anne}

\thanks{\textit{Date:} October 2008}

\begin{abstract}
We provide combinatorial models for all Kirillov--Reshetikhin crystals of nonexceptional type,
which were recently shown to exist.
For types $D_n^{(1)}$, $B_n^{(1)}$, $A_{2n-1}^{(2)}$ we rely on a previous construction 
using the Dynkin diagram automorphism which interchanges nodes $0$
and $1$. For type $C_n^{(1)}$  we use a Dynkin diagram folding and for types
$A_{2n}^{(2)}$, $D_{n+1}^{(2)}$ a similarity construction.
We also show that for types $C_n^{(1)}$ and $D_{n+1}^{(2)}$ the analog of the Dynkin diagram
automorphism exists on the level of crystals.
\end{abstract}

\maketitle

\section{Introduction}
Let $\geh$ be an affine Lie algebra and $U'_q(\geh)$ the corresponding quantum algebra without derivation.
Irreducible finite-dimensional $U_q'(\geh)$-modules were classified by Chari
and Pressley~\cite{CP,CP:1998} in terms of Drinfeld polynomials. It was then conjectured 
by Hatayama et al.~\cite{HKOTT:2002,HKOTY:1999} that a certain subset of such modules known 
as Kirillov--Reshetikhin (KR) modules $W^{(r)}_s$ have a crystal basis $B^{r,s}$ in the sense of 
Kashiwara \cite{Ka:1991}.
Here the index $r$ corresponds to a node of the Dynkin diagram of $\geh$ except the affine node
$0$ as specified in \cite{Kac}, and $s$ is an arbitrary positive integer. 
This conjecture was recently confirmed 
in~\cite{OS:2008} for all $\geh$ of nonexceptional affine type by using the results 
\cite{H:2006,H:2007,Nakajima:2003} on $T$-systems. (For many special cases including exceptional 
ones the conjecture was already known to be true 
in~\cite{BFKL,HN,JMO,KMN2:1992,Ka3,KMOY,Koga:1999,O:2007,Y:1998}.)
By the theory of affine finite crystals developed 
in~\cite{KMN1:1992,KMN2:1992}, that is, crystal bases of finite-dimensional $U'_q(\geh)$-modules, any 
integrable highest weight $U_q(\geh)$-module can be realized as a semi-infinite tensor
product of perfect crystals. This is known as the path realization.
Many of the crystals coming from KR modules, KR crystals for short, are (conjectured to be) 
perfect~\cite{HKOTT:2002,HKOTY:1999,KMN2:1992,S:2008}. By~\cite{FL:2007,NS:2006}
perfect KR crystals are isomorphic as classical crystals to certain Demazure
subcrystals of integrable highest weight crystals. In~\cite{FSS:2007} it was shown that under certain 
assumptions the classical isomorphism from the Demazure crystal to the KR crystal,
sends zero arrows to zero arrows. This implies in particular that the affine
crystal structure on these KR crystals is unique.

In this paper we solve the long outstanding problem of the construction of KR crystals.
We provide an explicit combinatorial crystal structure for all KR crystals $B^{r,s}$ of $W^{(r)}_s$
for $\geh$ of nonexceptional type. To do this we first
construct a combinatorial model $V^{r,s}(=V^{r,s}_\geh)$ for the KR crystal $B^{r,s}$. Let us look at 
type $A_{n-1}^{(1)}$ for instance. Since it is known that $W^{(r)}_s$ is irreducible as a 
$U_q(A_{n-1})$-module with highest weight $s\La_r$, the combinatorial crystal $V^{r,s}$ is defined 
to be the highest weight $A_{n-1}$-crystal $B(s\La_r)$, which can be identified with the set of 
semi-standard tableaux of rectangular shape $(s^r)$. On $B(s\La_r)$ the action of crystal operators 
$e_i,f_i$ ($i=1,2,\ldots,n-1$) is known \cite{KN:1994}. Hence, we are left to define the action of 
$e_0,f_0$. This was done by Shimozono~\cite{Sh:2002}, exploiting the fact that there is an
automorphism $\sigma$ defined on $B(s\La_r)$ which corresponds to the Dynkin diagram automorphism mapping $i$ to $i+1$ modulo $n$ (see Section~\ref{subsec:A}). With this $\sigma$
the affine crystal operator is given by $e_0=\sigma^{-1}\circ e_1 \circ \sigma$ and 
$f_0=\sigma^{-1}\circ f_1 \circ \sigma$.

\begin{table}
\begin{eqnarray*}
A_{n-1}^{(1)}
&\vcenter{\xymatrix@R=1ex{
&&*{\circ}<3pt> \ar@{-}[drr]^<{\;\,0} \ar@{-}[dll] \\
*{\circ}<3pt> \ar@{-}[r]_<{1} &*{\circ}<3pt> \ar@{-}[r]_<{2} 
&{} \ar@{.}[r]&{} \ar@{-}[r] &*{\circ}<3pt> \ar@{}[r]_<{n-1} &{}}}& (\phi,A_{n-1}) \\
B_n^{(1)}
&\vcenter{\xymatrix@R=1ex{
*{\circ}<3pt> \ar@{-}[dr]^<{0} \\
& *{\circ}<3pt> \ar@{-}[r]_<{2} 
& {} \ar@{.}[r]&{}  \ar@{-}[r]_>{\,\,\,\,n-1} &
*{\circ}<3pt> \ar@{=}[r] |-{\scalebox{2}{\object@{>}}}& *{\circ}<3pt>\ar@{}_<{n} \\
*{\circ}<3pt> \ar@{-}[ur]_<{1}}}& (\vdom,B_n) \\
C_n^{(1)}
&\vcenter{\xymatrix@R=1ex{
*{\circ}<3pt> \ar@{=}[r] |-{\scalebox{2}{\object@{>}}}_<{0} 
&*{\circ}<3pt> \ar@{-}[r]_<{1} 
& {} \ar@{.}[r]&{}  \ar@{-}[r]_>{\,\,\,\,n-1} &
*{\circ}<3pt> \ar@{=}[r] |-{\scalebox{2}{\object@{<}}}
& *{\bullet}<3pt>\ar@{}_<{n}}}& (\hdom,C_n) \\
D_n^{(1)}
&\vcenter{\xymatrix@R=1ex{
*{\circ}<3pt> \ar@{-}[dr]^<{0}&&&&&*{\bullet}<3pt> \ar@{-}[dl]^<{n-1}\\
& *{\circ}<3pt> \ar@{-}[r]_<{2} 
& {} \ar@{.}[r]&{} \ar@{-}[r]_>{\,\,\,n-2} &
*{\circ}<3pt> & \\
*{\circ}<3pt> \ar@{-}[ur]_<{1}&&&&&
*{\bullet}<3pt> \ar@{-}[ul]_<{n}}}& (\vdom,D_n) \\
A_{2n}^{(2)}
&\vcenter{\xymatrix@R=1ex{
*{\circ}<3pt> \ar@{=}[r] |-{\scalebox{2}{\object@{<}}}_<{0} 
&*{\circ}<3pt> \ar@{-}[r]_<{1} 
& {} \ar@{.}[r]&{}  \ar@{-}[r]_>{\,\,\,\,n-1} &
*{\circ}<3pt> \ar@{=}[r] |-{\scalebox{2}{\object@{<}}}
& *{\circ}<3pt>\ar@{}_<{n}}}& (\cell,C_n) \\
A_{2n-1}^{(2)}
&\vcenter{\xymatrix@R=1ex{
*{\circ}<3pt> \ar@{-}[dr]^<{0} \\
& *{\circ}<3pt> \ar@{-}[r]_<{2} 
& {} \ar@{.}[r]&{}  \ar@{-}[r]_>{\,\,\,\,n-1} &
*{\circ}<3pt> \ar@{=}[r] |-{\scalebox{2}{\object@{<}}}& *{\circ}<3pt>\ar@{}_<{n} \\
*{\circ}<3pt> \ar@{-}[ur]_<{1}}}& (\vdom,C_n) \\
D_{n+1}^{(2)}
&\vcenter{\xymatrix@R=1ex{
*{\circ}<3pt> \ar@{=}[r] |-{\scalebox{2}{\object@{<}}}_<{0} 
&*{\circ}<3pt> \ar@{-}[r]_<{1} 
& {} \ar@{.}[r]&{}  \ar@{-}[r]_>{\,\,\,\,n-1} &
*{\circ}<3pt> \ar@{=}[r] |-{\scalebox{2}{\object@{>}}}
& *{\bullet}<3pt>\ar@{}_<{n}}}& (\cell,B_n)
\end{eqnarray*}
\caption{\label{tab:Dynkin}Dynkin diagrams}
\end{table}

For other types of nonexceptional algebras, KR modules are not necessarily irreducible.
Let $\geh_0$ be the finite-dimensional simple Lie algebra obtained by removing $0$ from the Dynkin
diagram of $\geh$. 
In general $W^{(r)}_s$ decomposes into
\begin{equation} \label{eq:decomp W}
	W^{(r)}_s\simeq\bigoplus_\la V(\la)
\end{equation}
as a $U_q(\geh_0)$-module, where $V(\la)$ stands for the irreducible $U_q(\geh_0)$-module 
with highest weight corresponding to $\la$ and the sum runs over all partitions $\la$ that 
can be obtained from the $r\times s$ (or $r\times (s/2)$ only when $\geh=B_n^{(1)}$ and $r=n$) 
rectangle by removing pieces of shape $\nu$ (where $(\nu,\geh_0)$ are given in 
Table~\ref{tab:Dynkin}). There are some exceptions, in which case $W^{(r)}_s$ is irreducible
and does not decompose as in~\eqref{eq:decomp W}. We call these nodes $r$ ``exceptional'';
they are the filled nodes in Table~\ref{tab:Dynkin}.
\footnote{In~\cite{OS:2008} the node $r=n$ for type $B_n^{(1)}$ was accidentally
marked as exceptional.}
These decompositions were proven by Chari~\cite{Chari:2001} in the untwisted cases. In general
they can be proven from the results by Nakajima and Hernandez~\cite{H:2006,H:2007,Nakajima:2003}. 
See \cite{HKOTT:2002,HKOTY:1999}.

The combinatorial crystal $V^{r,s}$ is constructed according to whether $r$ is exceptional or not. 
We first treat the nonexceptional cases.
Since we already know the $\geh_0$-crystal structure via Kashiwara-Nakashima 
tableaux~\cite{KN:1994}, it is sufficient to define an appropriate automorphism related to a 
Dynkin diagram automorphism or an explicit $0$-action.
For types $D_n^{(1)}$, $B_n^{(1)}$, $A_{2n-1}^{(2)}$ we rely on the construction in~\cite{S:2008}
of an automorphism $\sigma$ which fixes the $\{2, 3,\ldots,n\}$-crystal structure and 
interchanges nodes $0$ and $1$, and define the affine crystal operator as 
$e_0 = \sigma \circ e_1 \circ \sigma$ (see Section~\ref{subsec:DBA}).
The $C_n^{(1)}$ crystal $V^{r,s}$ is realized as a ``virtual'' crystal inside the type $A_{2n+1}^{(2)}$ 
KR crystal $V^{r,s}$ using a folding of the Dynkin diagram (see Section~\ref{subsec:C}).
For $\geh$ of type $D_{n+1}^{(2)}$ or $A_{2n}^{(2)}$, Kashiwara's similarity method~\cite{Ka:1996}
is used to construct $V^{r,s}_\geh$ through a unique injective embedding
$S:V^{r,s}_\geh \to V^{r,s}_{C_n^{(1)}}$ (see Section~\ref{subsec:A(2)D(2)}).
The combinatorial crystals for the exceptional nodes such as $V^{n,s}$ for types $C_n^{(1)}$, 
$D_{n+1}^{(2)}$ and $V^{n,s}$, $V^{n-1,s}$ for type $D_n^{(1)}$ are treated in 
Section~\ref{sec:exceptional}.

The main theorem of this paper can be stated as follows:
\begin{theorem}
The combinatorial crystal $V^{r,s}$ given in this paper is isomorphic as a $\geh$-crystal to
the KR crystal $B^{r,s}$.
\end{theorem}
This theorem summarizes Theorem~\ref{thm:A} for type $A_{n-1}^{(1)}$ shown in \cite{KMN2:1992,Sh:2002},
Theorem~\ref{thm:DBA} for types $D_n^{(1)}$, $B_n^{(1)}$, $A_{2n-1}^{(2)}$ shown in \cite{OS:2008},
Theorem~\ref{thm:B} for $r=n$ for type $B_n^{(1)}$,
Theorem~\ref{thm:uniqueness} for types $C_n^{(1)}$, $D_{n+1}^{(2)}$, $A_{2n}^{(2)}$,
Theorems~\ref{thm:exceptional CD} for the exceptional node $r=n$ of types $C_n^{(1)}$, $D_{n+1}^{(2)}$,
and finally Theorem~\ref{thm:exceptional D} for exceptional nodes $r=n-1,n$ of type $D_n^{(1)}$.

The general strategy to deduce $V^{r,s} \cong B^{r,s}$ is to show a certain uniqueness theorem.
For example, for types $D_n^{(1)}$, $B_n^{(1)}$, $A_{2n-1}^{(2)}$ it is stated as follows. 
If $V^{r,s}$ and $B$ have the same decompositions as $\{1,2,\ldots,n\}$ and $\{0,2,\ldots,n\}$-crystals,
then they have to be isomorphic (see Section \ref{subsec:unique DBA}).

By construction, the Dynkin diagram automorphism $\sigma$ for type $A_{n-1}^{(1)}$, $B_n^{(1)}$, 
$D_n^{(1)}$, and $A_{2n-1}^{(2)}$ acts on the combinatorial crystal $V^{r,s}$.
The Dynkin diagrams for type $C_n^{(1)}$ and $D_{n+1}^{(2)}$ also have an automorphism
mapping $i\mapsto n-i$ for all $i\in \{0,1,\ldots,n\}$. However, from the construction of $V^{r,s}$
for these types using Dynkin diagram foldings and similarity methods, it is not obvious that
this Dynkin diagram automorphism extends to $V^{r,s}$. This is proven in Theorem~\ref{thm:auto}
and shows in particular that~\cite[Assumption 1]{OSS:2003} holds.

\subsection*{Organization}
The paper is organized as follows. In Section~\ref{sec:classical} we review some general facts and 
definitions about crystals, in particular the classical crystals of type $B_n$, $C_n$, $D_n$ using 
Kashiwara--Nakashima tableaux~\cite{KN:1994}.
In Section~\ref{sec:classical results} we review the branching $X_n\to X_{n-1}$ in terms of
$\pm$-diagrams, and derive some properties of $B_n$, $C_n$ crystals and their corresponding
$\pm$-diagrams. These definitions and properties are used to define the combinatorial KR
crystals $V^{r,s}$ in Section~\ref{sec:KR} and to show in Section~\ref{sec:uniqueness} that there is 
a unique crystal with the classical decompositions of $B^{r,s}$, thereby proving that 
$V^{r,s}\cong B^{r,s}$. The $V^{r,s}$ for exceptional nodes are treated in Section~\ref{sec:exceptional}.
In Section~\ref{sec:auto} it is shown that the Dynkin diagram automorphism of type $C_n^{(1)}$ and 
$D_{n+1}^{(2)}$ extends to $V^{r,s}$.

\subsection*{Acknowledgements.}
GF was supported in part by DARPA and AFOSR through the grant FA9550-07-1-0543 and by
the DFG-Projekt ``Kombinatorische Beschreibung von Macdonald und Kostka-Foulkes
Polynomen''.
MO was supported by grant JSPS 20540016.
AS was partially supported by the NSF grants DMS--0501101, DMS--0652641, and DMS--0652652.

GF and AS would like to thank the program ``Combinatorial representation theory''  held at MSRI
from January through May 2008, where part of this research was carried out.
MO and AS would like to thank the organizers of the conference
``Quantum affine Lie algebras, extended affine Lie algebras, and applications''
held at Banff where part of this work was carried out and presented. 
The implementation (by one of the authors) of crystals and in particular KR crystals in 
MuPAD-Combinat~\cite{HT:2003} and Sage~\cite{Sage} was extremely useful in undertaking 
the research for this article.

\section{Some review of crystal theory}
\label{sec:classical}

We review some basic definitions and facts about crystals that are used in this paper
in Section~\ref{subsec:general}. In order to describe the crystal graphs for the finite-dimensional modules 
of quantum groups of classical type, Kashiwara and Nakashima~\cite{KN:1994} introduced 
the analogue of semi-standard tableaux, called Kashiwara--Nakashima (KN) tableaux. In 
Sections~\ref{subsec:KN C}-\ref{subsec:KN D} we review KN tableaux for types $B_n$, $C_n$, and $D_n$,
respectively.

\subsection{General definitions}
\label{subsec:general}

Crystal theory was introduced by Kashiwara~\cite{Ka:1991} which provides a combinatorial way
to study the representation theory of quantum algebras $U_q(\geh)$. In this paper $\geh$ stands for
a simple Lie algebra or affine Kac--Moody Lie algebra with index set $I$ and $U_q(\geh)$ is the 
corresponding quantum algebra. 
Axiomatically, a $\geh$-crystal is a nonempty set $B$ together with maps
\begin{equation*}
\begin{split}
	e_i, f_i &: B \to B \cup \{\emptyset\} \qquad \text{for $i\in I$,}\\
	\wt &: B \to P,
\end{split}
\end{equation*}
where $P$ is the weight lattice associated to $\geh$. The maps $e_i$ and $f_i$ are Kashiwara's
crystal operators and $\wt$ is the weight function. Stembridge~\cite{St:2003} gave a local
characterization to determine when an axiomatic crystal actually corresponds to a 
$U_q(\geh)$-representation when $\geh$ is simply-laced. For further details about crystal theory, 
please consult for example~\cite{Ka:1991,HK:2002}.

To each crystal one can associate a crystal graph with vertices in $B$ and an arrow colored $i\in I$
from $b$ to $b'$ if $f_i(b)=b'$. For $b\in B$ and $i\in I$, let
\begin{equation*}
\begin{split}
	\ve_i(b) &= \max\{k \in \Z_{\ge 0} \mid e_i^k(b) \neq \emptyset \},\\
	\vp_i(b) &= \max\{k \in \Z_{\ge 0} \mid f_i^k(b) \neq \emptyset \}.
\end{split}
\end{equation*}
An element $b\in B$ is called highest (resp. lowest) weight if $e_i(b)=\emptyset$ (resp.
$f_i(b)=\emptyset$) for all $i\in I$.
For $J \subset I$, we say that $b\in B$ is $J$-highest (resp. $J$-lowest) if $e_i(b)=\emptyset$ (resp.
$f_i(b)=\emptyset$) for all $i\in J$.

We say that $b,b'\in B$ are $J$-related or $b\sim_J b'$ in symbols, if there exist
$J$-highest elements $b_0$, $b'_0$ of the same weight such that $b=f_{\vec{c}}(b_0)$,
$b'=f_{\vec{c}}(b'_0)$
for some sequence $\vec{c}$ from $J$. Here $f_{\vec{c}}=f_{c_1} \cdots f_{c_\ell}$ for
$\vec{c}=(c_1,\ldots,c_\ell)$.
A $J$-component $\mathcal{C}$ of a crystal $B$ is a connected component in the crystal graph of $B$
when only considering arrows colored $i\in J$.

We denote by $B(\La)$ the highest weight crystal of highest weight $\La$, where $\La$
is a dominant integral weight. Let $\La_i$ with $i\in I$ be the fundamental weights associated to a simple
Lie algebra of classical types, that is, $A_{n-1},B_n,C_n$ or $D_n$. 
Then as usual, a dominant integral weight $\La=\La_{i_1}+\cdots + \La_{i_k}$ is identified with a partition 
with columns of height $i_j$ for $1\le j\le k$, except when $\La_{i_j}$ is a spin weight, in which case 
we identify $\La_{i_j}$ with a column of height $n$ and width $1/2$. For type $B_n$ the fundamental
weight $\La_n$ is a spin weight and for type $D_n$ the fundamental weights $\La_{n-1}$ and
$\La_n$ are spin weights. In this paper we use French notation 
where parts are drawn in increasing order from top to bottom.
For type $A_{n-1}$, the highest weight crystal $B(\La)$ is given by the set of all semi-standard Young
tableaux of shape $\La$ over the alphabet $\{1,2,\ldots,n\}$. For types $B_n$, $C_n$, and $D_n$
elements in $B(\La)$ are given by Kashiwara--Nakashima (KN) tableaux~\cite{KN:1994}; they are
reviewed in the next subsections.

Let $B_1,B_2$ be crystals. Then $B_1\ot B_2=\{b_1\ot b_2\mid b_1\in B_1,b_2\in B_2\}$ can be endowed with
the structure of crystal. In order to compute the action of $e_i,f_i$ on multiple tensor products,
it is convenient to use the rule called ``signature rule". Let $b_1\ot b_2\ot\cdots\ot b_m$
be an element of the tensor product of crystals $B_1\ot B_2\ot\cdots\ot B_m$. One wishes to find 
the indices $j,j'$ such that 
\begin{align*}
e_i(b_1\ot\cdots\ot b_m)&=b_1\ot\cdots\ot e_ib_j\ot\cdots\ot b_m,\\
f_i(b_1\ot\cdots\ot b_m)&=b_1\ot\cdots\ot f_ib_{j'}\ot\cdots\ot b_m.
\end{align*}
To do it, we introduce ($i$-)signature by
\[
\overbrace{-\cdots-}^{\ve_i(b_1)}\overbrace{+\cdots+}^{\vp_i(b_1)}
\overbrace{-\cdots-}^{\ve_i(b_2)}\overbrace{+\cdots+}^{\vp_i(b_2)}
\:\cdots\cdots\:
\overbrace{-\cdots-}^{\ve_i(b_m)}\overbrace{+\cdots+}^{\vp_i(b_m)}.
\]
We then reduce the signature by deleting the adjacent $+-$ pair successively. 
Eventually we obtain a reduced signature of the following form.
\[
--\cdots-++\cdots+
\]
Then the action of $e_i$ (resp. $f_i$) corresponds to changing the rightmost $-$ to $+$ 
(resp. leftmost $+$ to $-$). If there is no $-$ (resp. $+$) in the signature, then the action of
$e_i$ (resp. $f_i$) should be set to $\emptyset$. The value of $\ve_i(b)$ (resp. $\vp_i(b)$)
is given by the number of $-$ (resp. $+$) in the reduced signature.

Consider, for instance, an element $b_1\ot b_2\ot b_3$ of the 3 fold tensor product $B_1\ot B_2\ot B_3$. 
Suppose $\ve_i(b_1)=1,\vp_i(b_1)=2,\ve_i(b_2)=1,\vp_i(b_2)=1,
\ve_i(b_3)=2,\vp_i(b_3)=1$. Then the signature and reduced one read 
\[
\begin{array}{cclcccr}
\mbox{sig}&&-++&\cdot&-+&\cdot&--+\phantom{.}\\
\mbox{red sig}&&-&\cdot&&\cdot&+.
\end{array}
\]
Thus we have
\begin{align*}
e_i(b_1\ot b_2\ot b_3)&=e_ib_1\ot b_2 \ot b_3,\\
f_i(b_1\ot b_2\ot b_3)&=b_1\ot b_2 \ot f_ib_3.
\end{align*}

\subsection{KN tableaux of type $C_n$}
\label{subsec:KN C}
In this section we review KN tableaux of type $C_n$.
On the set of letters $\{i,\overline{i} \mid 1\le i\le n\}$, introduce the
following order
\[
   1\prec2\prec\cdots\prec
n\prec\overline{n}\prec\cdots\prec\overline{2}\prec\overline{1}.
\]
As a set the crystal $B(\La_N)$ of the fundamental representation with highest weight $\La_N$
is given by
\begin{equation} \label{col cond C}
  \def\t{\tableau[sY]}
  \let\c=\vcenter
  B(\La_N)=\left\{
  \begin{array}{|c|} \hline i_N \\ \hline \vdots \\ \hline i_1 \\ \hline
\end{array} \mid
  \begin{array}{l}
   (1)\;1\preceq i_1\prec\cdots\prec i_N\preceq\ol{1},\\
   (2)\;\text{if $i_k=p$ and $i_l=\overline{p}$, then $k+(N-l+1)\le p$}
   \end{array}
   \right\}.
\end{equation}
To describe the crystal $B(\La_M+\La_N)$ ($M\ge N$) we need to define the notion
of $(a,b)$-configurations. 

\begin{definition} \label{ab-config C}
Let
\[
  u=\begin{array}{|c|} \hline i_M \\ \hline \vdots \\ \hline i_1 \\ \hline
\end{array}  \in B(\La_M)
  \quad\text{and}\quad
   v=     \begin{array}{|c|} \hline j_N \\ \hline \vdots \\ \hline j_1 \\ \hline
\end{array}  \in B(\La_N).
\]
For $1\le a\le b\le n$, we say $w=(u,v)$ is in the
$(a,b)$-configuration
if it satisfies the following: There exist $1\le p\le q<r\le s\le N$ such that
$i_p=a,i_q=b,
i_r=\overline{b},j_s=\overline{a}$ or
$i_p=a,j_q=b,j_r=\overline{b},j_s=\overline{a}$. The definition
includes the case where $a=b$, $p=q$, and $r=s$. Define
\[
   p(a,b;w)=(q-p)+(s-r).
\]
\end{definition}
Then the crystal $B(\La_M+\La_N)$ of the highest weight module of highest weight
$\La_M+\La_N$
is given by
\begin{equation} \label{adj cond C}
B(\La_M+\La_N) = \left\{
   w= \begin{array}{|c|c|} \cline{1-1} i_M & \multicolumn{1}{c}{} \\ \hline
\vdots & j_N \\ \hline
   \vdots & \vdots \\ \hline i_1 & j_1\\ \hline \end{array} \mid
   \begin{array}{l}
   (1)\;i_k\preceq j_k\text{ for }1\le k\le N,\\
   (2)\;\text{if $w$ is in the $(a,b)$-configuration,}\\
   \phantom{(2)}\; \text{then $p(a,b;w)<b-a$}
\end{array}
\right\}.
\end{equation}
Note that an element of $B(\La_M+\La_N)$ cannot be in the $(a,a)$-configuration.
We can now
describe the crystal $B(\La)$ of the highest weight module of highest weight
$\La=\La_{l_1}+\cdots+\La_{l_p}$ ($n\ge l_1\ge\cdots\ge l_p\ge1$) as

\begin{equation} \label{B(La)}
  B(\La) = \left\{
  w= \begin{array}{|c|c|c|c|} \cline{1-1}  & \multicolumn{3}{c}{} \\ 
\cline{2-2} & & \multicolumn{2}{c}{} \\ \cline{3-3}& \vdots& \vdots&
\multicolumn{1}{c}{} \\ \cline{4-4} t_1 & & &t_p \\  &  &  &\\ \cline{1-4}
\end{array} \mid
  \begin{array}{|c|c|} \cline{1-1}  & \multicolumn{1}{c}{} \\ \cline{2-2}  &  \\
t_k & t_{k+1} \\ & \\ \cline{1-2}  \end{array} \in
B(\La_{l_k}+\La_{l_{k+1}})\text{ for any }k=1,\ldots,p-1
\right\}.
\end{equation}

Let us describe the action of crystal operators on $B(\La)$. For the simplest case $B(\La_1)$
they are given by the following crystal graph.
\[
\tab{1}\overset{1}{\longrightarrow}\tab{2}\overset{2}{\longrightarrow}\cdots
\overset{n-1}{\longrightarrow}\tab{n}\overset{n}{\longrightarrow}\tab{\ol{n}}
\overset{n-1}{\longrightarrow}\cdots\overset{2}{\longrightarrow}\tab{\ol{2}}
\overset{1}{\longrightarrow}\tab{\ol{1}}
\]
For the general case, we regard a tableau as an element of $B(\La_1)^{\ot N}$, where $N$ is the 
number of boxes of the tableau. We move along the tableau from the rightmost column 
to left, and in each column we move from bottom to top. Then we obtain the sequence of letters
$b_1,b_2,\ldots,b_N$. We associate the tableau to $b_1\ot b_2\ot\cdots\ot b_N$ in 
$B(\La_1)^{\ot N}$. Then the action of $e_i,f_i$ is given by the multiple tensor product rule
explained in the previous subsection.

\subsection{KN tableaux of type $B_n$}
\label{subsec:KN B}
The construction for $B_n$ is divided into two cases.
Set $\omega_i = \Lambda_i$ (for $i \leq n-1$), $\omega_n = 2 \Lambda_n$, and on the
set of letters $\{i,\overline{i} \mid 1\le i\le n\} \cup \{0\}$ introduce the following order
\[
	1\prec2\prec\cdots\prec n\prec 0 \prec
\overline{n}\prec\cdots\prec\overline{2}\prec\overline{1}.
\]
As a set the crystal $B( \omega_N)$ of the fundamental representation with
highest weight $\omega_N$ is given by
\begin{equation} \label{col cond B}
  \def\t{\tableau[sY]}
  \let\c=\vcenter
  B(\omega_N)=\left\{
  \begin{array}{|c|} \hline i_N \\ \hline \vdots \\ \hline i_1 \\ \hline
\end{array} \mid
  \begin{array}{l}
  (1)\;1\preceq i_1\prec \cdots \prec i_N\preceq\ol{1},\\
  \text{ but no element other than $0$ can appear more than once}\\
  (2)\;\text{if $i_k=p$ and $i_l=\overline{p}$, then $k+(N-l+1)\le p$}
  \end{array}
  \right\}.
\end{equation}
The second case is the ''spin representation'' of highest weight $\Lambda_n$.
Define on $\{1, \ldots, n,\ol{n}, \ldots, \ol{1}\}$ a linear order by
$$1\prec2\prec\cdots\prec n\prec
\overline{n}\prec\cdots\prec\overline{2}\prec\overline{1}.$$
Then as a set the crystal $B(\Lambda_n)$ is given by
\begin{equation} \label{col cond B_spin}
  \def\t{\tableau[sY]}
  \let\c=\vcenter
  B(\La_n)=\left\{
  \begin{array}{|c|} \hline i_n \\ \hline \vdots \\ \hline i_1 \\ \hline
\end{array} \mid
  \begin{array}{l}
  (1)\;1\preceq i_1\prec\cdots\prec i_n\preceq\ol{1},\\
  (2) \text{ $i$ and $\ol{i}$ do not appear simultaneously }\\
  \end{array}
  \right\}.
\end{equation}

To describe the crystal $B(\omega_M + \omega_N)$ ($M\ge N$) we again need to define
the notion of $(a,b)$-configurations. 
\begin{definition}\label{ab-config B}
Let
\[
  u=     \begin{array}{|c|} \hline i_M \\ \hline \vdots \\ \hline i_1 \\ \hline
\end{array}  \in B(\omega_M)
  \quad\text{and}\quad
  v=     \begin{array}{|c|} \hline j_N \\ \hline \vdots \\ \hline j_1 \\ \hline
\end{array}  \in B(\omega_N).
\]
For $1\le a \le b < n$ we have the same conditions as in Definition~\ref{ab-config C} for type $C_n$.
For $1\le a < n$, we say $w=(u,v)$ is in the $(a,n)$-configuration
if it satisfies the following: There exist $1\le p\le q<r = q+1 \le s \le N$
such that $i_p=a$, $j_s=\overline{a}$ and one of the conditions is satisfied:
\begin{enumerate}
\item $i_q$ and $i_r ( = i_{q+1})$ are $n, 0, $ or $\ol{n}$.
\item $j_q$ and $j_r (= j_{q+1})$ are $n,0, $ or $\ol{n}$.
\end{enumerate}
We say $w=(u,v)$ is in the $(n,n)$-configuration
if there are $1\le p <  q \le N$ such that $i_p=n $ or $0$ and $j_p = 0$ or
$\ol{n}$.
Define again
\[
  p(a,b;w)=(q-p)+(s-r).
\]
\end{definition}
Then as a set the crystal $B(\omega_M + \omega_N)$ of the highest weight module of
highest weight $\omega_M + \omega_N$ is given by
\begin{equation} \label{adj cond B}
  B(\omega_M + \omega_N) = \left\{
  w= \begin{array}{|c|c|} \cline{1-1} i_M & \multicolumn{1}{c}{} \\ \hline
\vdots & j_N \\ \hline
  \vdots & \vdots \\ \hline i_1 & j_1\\ \hline \end{array} \mid
  \begin{array}{l}
  (1)\;i_k\preceq j_k\text{ for }1\le k\le N,\\
  \phantom{(1)}\; \text{ and $i_k$ and $j_k$ cannot both be 0,}\\
  (2)\;\text{if $w$ is in the $(a,b)$-configuration,}\\
  \phantom{(2)}\; \text{then $p(a,b;w)<b-a$}
  \end{array}
\right\}.
\end{equation}
Note that an element of $B(\omega_M + \omega_N)$ cannot be in the
$(a,a)$-configuration.
The conditions for $B(\Lambda_n + \omega_N)$ are formally just the same as in \eqref{adj cond B},
although there is no $0$ in the first column.

We can now describe the crystal $B(\La)$ of the highest weight module of highest weight
$\La$. If $\La$ is of the form $\La=\omega_{l_1}+\cdots+\omega_{l_p}$ ($n \ge l_1\ge\cdots\ge l_p\ge1$), 
\begin{equation}
  B(\La) = \left\{
  w= \begin{array}{|c|c|c|c|} \cline{1-1}  & \multicolumn{3}{c}{} \\ 
\cline{2-2} & & \multicolumn{2}{c}{} \\ \cline{3-3}& \vdots& \vdots&
\multicolumn{1}{c}{} \\ \cline{4-4} t_1 & & &t_p \\  &  &  &\\ \cline{1-4}
\end{array} \mid
  \begin{array}{|c|c|} \cline{1-1}  & \multicolumn{1}{c}{} \\ \cline{2-2}  &  \\
t_k & t_{k+1} \\ & \\ \cline{1-2}  \end{array} \in  B(\omega_{l_k} +
\omega_{l_{k+1}})\text{ for }k \le p-1
\right\}.
\end{equation}
The crystal structure can be described in the same way as in Section \ref{subsec:KN C}. The only
difference is that we have to replace the crystal graph of $B(\La_1)$ with the one below.
\[
\tab{1}\overset{1}{\longrightarrow}\tab{2}\overset{2}{\longrightarrow}\cdots
\overset{n-1}{\longrightarrow}\tab{n}\overset{n}{\longrightarrow}
\tab{0}\overset{n}{\longrightarrow}\tab{\ol{n}}\overset{n-1}{\longrightarrow}\cdots
\overset{2}{\longrightarrow}\tab{\ol{2}}\overset{1}{\longrightarrow}\tab{\ol{1}}
\]

Else, $\La$ can be written as $\La=\La_n+\omega_{l_2}+\cdots+\omega_{l_p}$ ($n \ge l_2\ge\cdots\ge l_p\ge1$).
In this case, 
\begin{equation}
  B(\La) = \left\{
  w= \begin{array}{|c|c|c|c|} \cline{1-1}  & \multicolumn{3}{c}{} \\ 
\cline{2-2} & & \multicolumn{2}{c}{} \\ \cline{3-3}& \vdots& \vdots&
\multicolumn{1}{c}{} \\ \cline{4-4} t_1 & & &t_p \\  &  &  &\\ \cline{1-4}
\end{array} \mid
  \begin{array}{|c|c|} \cline{1-1}  & \multicolumn{1}{c}{} \\ \cline{2-2}  &  \\
t_k & t_{k+1} \\ & \\ \cline{1-2}  \end{array} \in \begin{array}{l}
B(\Lambda_n+\omega_{l_2})\text{ for }k=1 \\
B(\omega_{l_k} + \omega_{l_{k+1}})\text{ for }2\le k\le p-1
\end{array}
\right\}.
\end{equation}
The crystal structure in this case is given as that on the tensor product 
$B(\omega_{l_2}+\cdots+\omega_{l_p})\ot B(\La_n)$. For the crystal structure on $B(\La_n)$ 
see \cite[Section 5.4]{KN:1994}.

\subsection{KN tableaux of type $D_n$}
\label{subsec:KN D}
For type $D_n$ we only consider $B(\La)$, where the coefficients of $\La_{n-1}$ and $\La_n$
of $\La$ are $0$. On the
set $\{i ,\ol{i} \mid \; 1 \leq i \leq n \}$ we define the following order
$$ 1 \prec 2 \prec \ldots \prec n-1 \prec {n \atop \ol{n}} \prec \ol{n-1} \prec
\ldots \prec \ol{1} $$ where there is no order between $n$ and $\ol{n}$.
Then as a set the crystal of highest weight $\La_N$ for  $N \leq n-2$ is
\begin{equation} \label{col cond D1}
  \def\t{\tableau[sY]}
  \let\c=\vcenter
  B(\La_N)=\left\{
  \begin{array}{|c|} \hline i_N \\ \hline \vdots \\ \hline i_1 \\ \hline
\end{array} \mid
  \begin{array}{l}
  (1)\;i_{j}\not\succeq i_{j+1} \text{ for } 1 \leq j < N\\
  (2)\;\text{if $i_k=p$ and $i_l=\overline{p}$ ($1 \leq p \leq n$), then $k+(N-l+1)\le p$}
  \end{array}
  \right\}.
\end{equation}

To describe the crystal $B(\La_M + \La_N)$ ($M\ge N$) we need to define
again the notion of $(a,b)$-configurations. 

\begin{definition}\label{ab-config D}\mbox{}
Let
\[
  u=     \begin{array}{|c|} \hline i_M \\ \hline \vdots \\ \hline i_1 \\ \hline
\end{array}  \in B(\La_M)
  \quad\text{and}\quad
  v=     \begin{array}{|c|} \hline j_N \\ \hline \vdots \\ \hline j_1 \\ \hline
\end{array}  \in B(\La_N).
\]
\begin{enumerate}
\item For $1\le a \le b<n$ we have the same conditions as in Definition~\ref{ab-config C}.
\item  For $1\le a  < n$, we say $w=(u,v)$ is in the $(a,n)$-configuration
if it satisfies the following: There exist $1\le p\le q<r = q+1 \le s \le N$
such that $i_p=a,j_s=\overline{a}$ and one of the conditions is satisfied:
\begin{enumerate}
\item $i_q$ and $i_r ( = i_{q+1})$ are $n$ or $\ol{n}$,
\item $j_q$ and $j_r (= j_{q+1})$ are $n $ or $\ol{n}$.
\end{enumerate}
\item We say $w=(u,v)$ is in the $(n,n)$-configuration
if there are $1\le p <  q \le N$ such that $i_p=n $ or $\ol{n}$ and $j_p = n$ or $\ol{n}$.
\item For $1 \leq a < n$, $w=(u,v)$ is in the $a$-odd-configuration if the following conditions
are satisfied: There exists $1 \leq p \leq q < r \leq s \leq N$ such that
\begin{enumerate}
\item $r-q +1 $ is odd,
\item $i_p = a$ and $j_s = \ol{a}$,
\item $j_q = n, i_r = \ol{n}$ or $ j_q = \ol{n}, i_r = n$.
\end{enumerate}
\item For $1 \leq a < n$, $w=(u,v)$ is in the $a$-even-configuration
if the following conditions are satisfied: There exists $1 \leq p \leq q < r \leq s \leq N$ such that
\begin{enumerate}
\item $r-q +1 $ is even,
\item $i_p = a$ and $j_s = \ol{a}$,
\item $j_q = n, i_r = n$ or $ j_q = \ol{n}, i_r = \ol{n}$.
\end{enumerate}
\end{enumerate}
Then
\begin{enumerate}
\item If $w$ is in the $(a,b)$-configuration for $1 \leq a \leq b \leq n$, we
define $p(a,b;w) = (q-p) + (s-r)$. If $a=b=n$, set $p(a,b;w) = 0$.
\item If $w$ is in the $a$-odd or $a$-even-configuration, we define $q(a;w) = s-p$.
\end{enumerate}
\end{definition}
Then the crystal $B(\La_M + \La_N)$ of the highest weight module of
highest weight $\La_M + \La_N$ ($n\ge M\ge N\ge 1$) is given by
\begin{equation} \label{adj cond D}
  B(\La_M + \La_N) = \left\{
  w= \begin{array}{|c|c|} \cline{1-1} i_M & \multicolumn{1}{c}{} \\ \hline
\vdots & j_N \\ \hline
  \vdots & \vdots \\ \hline i_1 & j_1\\ \hline \end{array} \mid
  \begin{array}{l}
  (1)\;i_k\preceq j_k\text{ for }1\le k\le N,\\
  (2)\;\text{if $w$ is in the $(a,b)$-configuration,}\\
  \phantom{(2)}\; \text{then $p(a,b;w)<b-a$}\\
  (3)\;\text{if $w$ is in the $a$-odd-} \\
  \phantom{(3)}\; \text{or $a$-even-configuration,}\\
  \phantom{(3)}\; \text{then $q(a;w) < n-a$}
  \end{array}
\right\}.
\end{equation}
The crystal $B(\La)$ with $\La=\La_{l_1}+\cdots+\La_{l_p}$ ($n-2\ge l_1\ge l_2\ge\cdots\ge l_p\ge1$)
is described again as \eqref{B(La)}. The crystal structure on $B(\La)$ is obtained as in Section 
\ref{subsec:KN C}. The only difference is that we have to replace the crystal graph of $B(\La_1)$ 
with the one below.
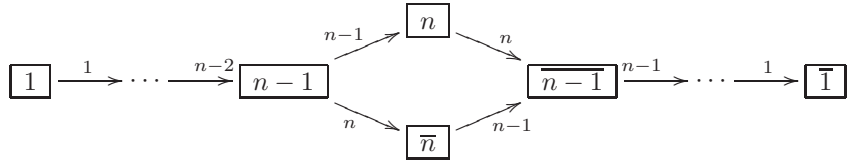
\[
\vcenter{\xymatrix@R=1ex{
&&&\hbox{$\tab{n}$}\ar[dr]^n\\
\hbox{$\tab{1}$}\ar[r]^1&\cdots\ar[r]^{n-2}&\hbox{$\tab{n-1}$}\ar[ur]^{n-1}\ar[dr]_n&&
\hbox{$\tab{\ol{n-1}}$}\ar[r]^{n-1}&\cdots\ar[r]^1&\hbox{$\tab{\ol{1}}$}\\
&&&\hbox{$\tab{\ol{n}}$}\ar[ur]_{n-1}}}
\]

\section{Properties and branching of classical crystals}
\label{sec:classical results}

In this section we derive some properties of $B_n$, $C_n$ crystals in Section~\ref{subsec:prop BC},
branching rules $X_n\to X_{n-1}$ in terms of $\pm$-diagrams where $X=B,C,D$ in 
Section~\ref{subsec:branching}, and properties of $\pm$-diagrams  for type $B_n$, $C_n$ in 
Section~\ref{subsec:pm properties BC}. These results will be used later in the construction of KR
crystals.

\subsection{Properties of $B_n$ and $C_n$ crystals}
\label{subsec:prop BC}
In this section we prove some preliminary results for the form and properties of special elements
in $B(\La)$ of type $C_n$ and $B_n$. In addition to the index set of the classical Dynkin diagram
$I_0=\{1,2,\ldots,n\}$, we will also use $J=\{1,2,\ldots,n-1\}$ and $J'=\{2,3,\ldots,n\}$.

\begin{lemma} \label{lemma:Jlowest}
Let $b$ be a $J$-lowest weight element of type $C_n$ or $B_n$.
Then the columns of $b$ must be of the form
\begin{equation} \label{eq:Jlowest}
\begin{array}{|c|}
\hline
\ol{l}_t^{\phantom{L}}\\ \hline
\vdots\\ \hline 
\ol{l}_1^{\phantom{L}}\\ \hline
0^\alpha\\ \hline
n\\ \hline
n-1\\ \hline
\vdots\\ \hline
k \\ \hline
\end{array}
\end{equation}
where $l_i<k$ for all $1\le i\le t$, $\alpha\ge 0$, and $\alpha=0$ for type $C_n$.
\end{lemma}

\begin{proof}
We prove this lemma for $C_n$. The $B_n$ case can be treated in exactly the same way.

We prove the lemma by induction on the columns (from left to right).

All unbarred letters of $b$ form a tableau in the bottom left (in French convention).
If there is no unbarred letter in the first column, then there is none in the
entire tableau. Suppose the letter $k$ occurs in the first column, then there has to be (since $b$ is
lowest weight) a $k+1$ or a $\ol{k}$ in the same column. Suppose there is a $\ol{k}$. Then,
since $b$ is lowest weight, there has to be a $\ol{k-1}$. Then
the top of the column has the form
\begin{equation*}
\begin{array}{|c|}
\hline
\ol{1}^{\hspace{-2mm}\phantom{L}}\\ \hline
\ol{2}^{\hspace{-2mm}\phantom{L}}\\ \hline
\vdots\\ \hline
\ol{k}^{\hspace{-2mm}\phantom{L}} \\ \hline
\vdots \\ \hline
\end{array} \; .
\end{equation*}
But this is not a legal type $C_n$ tableau, since the pair $k$ and $\ol{k}$ does not satisfy the
condition of~\eqref{col cond C}.
Hence there must be a $k+1$. Therefore the first column looks as follows
\begin{equation*}
\begin{array}{|c|}
\hline
\ol{l}_r^{\hspace{-2mm}\phantom{L}}\\ \hline
\vdots\\ \hline
\ol{l}_1^{\hspace{-2mm}\phantom{L}}\\ \hline
n\\ \hline
\vdots\\ \hline
k+1\\ \hline
k \\ \hline
\end{array} \; .
\end{equation*}
The condition $l_i < k$ for all $i$ is forced again by the tableau rules~\eqref{col cond C} for $C_n$,
since the top has to be consecutive in the barred letters for $b$ to be $J$-lowest weight. This
concludes the induction beginning. Now suppose the claim is true for the first $m-1$ columns.
If the $m$-th column does not start with an unbarred letter, we are done.
So suppose it starts with a $k$.

Case 1: The column $m-1$ starts with a $k$ as well. Since $b$ is $J$-lowest weight,
there must be either a $k+1$ or a $\ol{k}$ to bracket the $k$ in column $m$. By induction, in
every column to the left of the $m$-th column, there is no $\ol{k}$. Suppose there is a $\ol{k}$ in
the $m$-th column. Then $b$ contains the following pattern:
\begin{equation} \label{eq:pattern}
\begin{array}{|c|c|}
\hline
* & \ol{k}^{\hspace{-2mm}\phantom{L}}\\ \hline
\vdots & \vdots \\ \hline
k &  k \\ \hline
\end{array}
\end{equation}
which is illegal by~\eqref{adj cond C}. Hence there must be a $k+1$. This letter $k+1$ has to
be bracketed again, since $b$ is lowest weight. Suppose it is bracketed with $\ol{k+1}$. But
by induction to the left there is no unbracketed $k+1$ or $\ol{k}$, which is a contradiction.
Hence there is a $k+2$ and so on.

Case 2: There is no column to the left of the $m$-th column, that starts with a $k$.
Suppose $k$ is bracketed with a $\ol{k}$. This $\ol{k}$ has to be in the same column as $k$,
since by induction, there is no $\ol{k}$ to the left. But this $\ol{k}$ has to be bracketed
by a $\ol{k-1}$ which is in the same column (by induction). Hence this column is of the form
\begin{equation*}
\begin{array}{|c|}
\hline
\ol{k-i}^{\hspace{-2mm}\phantom{L}}\\ \hline
\vdots\\ \hline
\ol{k-1}^{\hspace{-2mm}\phantom{L}}\\ \hline
\ol{k}^{\hspace{-2mm}\phantom{L}}\\ \hline
\vdots\\ \hline
k \\ \hline
\end{array}
\end{equation*}
for some $i$. The letter $\ol{k-i}$ has to be bracketed as well. Suppose it is bracketed with a $k-i$.
Then we obtain an illegal $C_n$ tableau by~\eqref{adj cond C}.
If it is bracketed with a $\ol{k-i-1}$ in another column, then there is a
corresponding unbarred letter in this column. This letter has to be bigger
or equal to $k-i$, which gives again an illegal tableau. The column
\begin{equation*}
\begin{array}{|c|}
\hline
\ol{1}^{\hspace{-2mm}\phantom{L}}\\ \hline
\vdots\\ \hline
\ol{k-1}^{\hspace{-2mm}\phantom{L}}\\ \hline
\ol{k}^{\hspace{-2mm}\phantom{L}}\\ \hline
\vdots\\ \hline
k \\ \hline
\end{array}
\end{equation*}
is also illegal by~\eqref{col cond C}. Hence $k$ has to be bracketed by $k+1$.
Repeating the arguments as before one obtains that the $m$-th column is of the form
\begin{equation*}
\begin{array}{|c|}
\hline
\vdots\\ \hline
n\\ \hline
\vdots\\ \hline
k+1\\ \hline
k \\ \hline
\end{array}
\end{equation*}

To finish the induction, we need to show that above the $n$ there are only $\ol{l}_i$
with $l_i < k$.
Suppose there is an $\ol{l}$ with $l\ge k$. Then, since the element is lowest weight, $\ol{l}$ must be
bracketed. It cannot be bracketed with an unbarred letter to the left of the $m$-th column,
since they are all smaller than $k$ (if there is an unbarred letter $l = k$, then the tableau
has a pattern of the form~\eqref{eq:pattern}, which is illegal by~\eqref{adj cond C}).
By induction the only barred letters to the left of the $m$-th column are $\ol{p}$, with $p < k-1$.
(If there is a $k$ at the left of the $m$-th column, then there could not be
a $\ol{k-1}$, because the tableau would be again illegal).
Hence all barred letters in the $m$-th column must be smaller than $k$.
\end{proof}

The inner shape of a $J$-lowest weight element $b$ is defined to be the shape after
deleting all $\ol{1}$s.

\begin{lemma} \label{lemma:inner shape}
The highest weight of the $J'$-component of a $J$-lowest weight element $b$
is given by the inner shape of $b$.
\end{lemma}

\begin{proof}
By Lemma \ref{lemma:Jlowest}, $b$ contains no $1$'s.
Construct a tableau $b'$ from $b$ by deleting all $\ol{1}$'s and replacing each letter $c$
(resp. $\ol{c}$) by $c-1$ (resp. $\ol{c-1}$). Then one finds that
$b'$ is again a KN tableau for $C_{n-1}$ ($B_{n-1}$ resp.): To be a KN tableau there are
conditions in each column~\eqref{col cond C} (resp.~\eqref{col cond B}) and those for adjacent
columns~\eqref{adj cond C} (resp.~\eqref{adj cond B}). The former condition is satisfied,
since there is no $(k,\ol{k})$ pair by Lemma~\ref{lemma:Jlowest}. The latter condition is invariant 
under changing contents as above.
So by applying a sequence of $e_i$ with $i\in\{1,2,\ldots,n-1\}$, $b'$ can be raised to the 
$C_{n-1}$-highest weight element ($B_{n-1}$ resp.).
\end{proof}

Recall that $b,b'\in B(\Lambda)$ are $J'$-related or $b\sim_{J'} b'$ in symbols, if there exist
$J'$-highest elements $b_0,b'_0$ of the same weight such that $b=f_{\vec{c}}(b_0)$,
$b'=f_{\vec{c}}(b'_0)$ for some sequence $\vec{c}$ from $J'$.

\begin{corollary} \label{cor:J-lt}
If $b_1$ and $b_2$ are $J$-lowest, $b_1 \sim_{J'} b_2$ and $\wt(b_1)=\wt(b_2)$, then $b_1$
and $b_2$ differ just by boxes containing $\ol{1}$. The inner tableaux are the same.
\end{corollary}
\begin{proof}
This follows directly from Lemma~\ref{lemma:inner shape}.
\end{proof}

\subsection{$X_n\to X_{n-1}$ branching and $\pm$-diagrams}
\label{subsec:branching}
Let $\geh_0$ be a finite Lie algebra of type $X_n=D_n,B_n$, or $C_n$.
In this section we describe a branching rule for $X_n\to X_{n-1}$ involving $\pm$-diagrams. 

A $\pm$-diagram $P$ of shape $\La/\la$ is a sequence of partitions $\la\subset \mu \subset \La$ 
such that $\La/\mu$ and $\mu/\la$ are horizontal strips (i.e. every column contains at most one box). We
depict this $\pm$-diagram by the skew tableau of shape $\La/\la$ in
which the cells of $\mu/\la$ are filled with the symbol $+$ and
those of $\La/\mu$ are filled with the symbol $-$. Write
$\La=\os(P)$ and $\la=\is(P)$ for the outer and inner shapes of the
$\pm$-diagram $P$. 
When drawing partitions or tableaux, we use the French
convention where the parts are drawn in increasing order from top to bottom.

There are a couple further type-specific requirements:
\begin{enumerate}
\item
For type $C_n$ the outer shape $\La$ contains columns of height at most $n$, but the inner 
shape $\lambda$ is not allowed to be of height $n$ (hence there are no empty columns of height $n$).
\item
For type $B_n$ the outer shape $\La$ contains columns of height at most $n$; for the columns
of height $n$, the $\pm$-diagram can contain at most one $0$ between $+$ and $-$ at height $n$
and no empty columns are allowed; furthermore there may be a spin column of height $n$ and 
width $1/2$ containing $+$ or $-$.
\item
For type $D_n$ suppose $\La= k_1 \La_1 + \cdots + k_{n-1} \La_{n-1} + k_n \La_n$.
If $k_n\ge k_{n-1}$ we depict this weight by $(k_n-k_{n-1})/2$ columns of height $n$ colored 1
(where we interpret a $1/2$ column as a $\La_n$ spin column if $k_n-k_{n-1}$ is odd),
$k_{n-1}$ columns of height $n-1$, and as usual $k_i$ columns of height $i$ for $1\le i\le n-2$.
If $k_n< k_{n-1}$ we depict this weight by $(k_{n-1}-k_n)/2$ columns of height $n$ colored 2
(where we interpret a $1/2$ column as a $\La_{n-1}$ spin column if $k_{n-1}-k_n$ is odd),
$k_n$ columns of height $n-1$, and as usual $k_i$ columns of height $i$ for $1\le i\le n-2$.
We require that columns of height $n$ are colored, contain $+$, $-$, or $\mp$,
but cannot simultaneously contain $+$ and $-$; spin columns can only contain $+$ or $-$.
\end{enumerate}

\begin{proposition} \label{P:branch}
For an $X_n$ dominant weight $\La$, there is an isomorphism of $X_{n-1}$-crystals
\begin{align*}
  B_{X_n}(\La) \cong \bigoplus_{\substack{\text{$\pm$-diagrams $P$} \\ \os(P)=\La}}
  B_{X_{n-1}}(\is(P)).
\end{align*}
That is, the multiplicity of $B_{X_{n-1}}(\la)$ in $B_{X_n}(\La)$,
is the number of $\pm$-diagrams of shape $\La/\la$.
\end{proposition}

\begin{proof}
This follows directly from the branching rules for $X_n$ to $X_{n-1}$
(see for example~\cite[pg. 426]{FH:1991}).
\end{proof}

There is a bijection $\bij:P\mapsto b$ from $\pm$-diagrams $P$ of shape 
$\La/\la$ to the set of $X_{n-1}$-highest weight vectors $b$ of $X_{n-1}$-weight $\la$.  
Namely, we construct a string of operators 
$f_{\aar} := f_{a_1} f_{a_2} \cdots f_{a_\ell}$ such that $\bij(P) = f_{\aar} u$, where $u$ is 
the highest weight vector in $B_{X_n}(\Lambda)$. Start with $\aar=()$. 
\begin{enumerate}
\item
Scan the columns of $P$ from right to left. For each column of $P$ for which a $+$ can be added, 
append $(1,2, \ldots, h)$ to $\aar$, where $h$ is the height of the added $+$. 
Note that a $+$ is not addable to a spin column. Further type specific rules are:
\begin{enumerate}
\item
For $D_n$, if there is an addable $+$ at height $n-1$ to a column of height $n$, append 
$(1,2,\ldots,n-2,n)$ if the color is $1$ and $(1,2,\ldots,n-1)$ if the color is $2$.
\item 
For $B_n$, if there is a column of height $n$ containing $0$ append $(1,2,\ldots,n)$ to $\aar$.
\end{enumerate}
\item
Next scan $P$ from left to right for columns containing a $-$ at height $h$.
\begin{enumerate}
\item
For $D_n$, if $h=n$, append $(1,2,\ldots,n-2,n)$ if the color is $1$ and 
$(1,2,\ldots, n-1)$ if the color is $2$ (this also applies to spin columns). If $h=n-1$, append 
$(1,2,\ldots,n)$. For $h<n-1$, append the string $(1,2,\ldots,n,n-2,n-3,\ldots, h)$ to $\aar$.
\item
For $B_n$, if the $-$ is in a spin column, append $(1,2, \ldots,n)$ to $\aar$.
Otherwise append the string $(1,2,\ldots,n-1,n,n,n-1,\ldots,h)$ to $\aar$.
\item
For $C_n$, append the string $(1,2,\ldots,n-1,n,n-1,\ldots,h)$ to $\aar$.
\end{enumerate}
\end{enumerate}
This correspondence can easily be checked by explicitly writing down the $X_{n-1}$-highest
weight vectors for each type.

Alternatively, suppose $\La$ is a dominant weight; we require that $\La$ does not contain any 
columns of height $n$ for type $D_n$.
Then the bijection $\bij:P\mapsto b$ from $\pm$-diagrams $P$ of shape 
$\La/\la$ to the set of $X_{n-1}$-highest weight vectors $b$ of
$X_{n-1}$-weight $\la$ is as follows. For any columns of height $n$ containing $+$, place a column
$12\ldots n$ (this includes spin columns for type $B_n$). Otherwise, place $\ol{1}$ in all positions in 
$P$ that contain a $-$, place a $0$ in the position containing $0$, and fill the remainder of all 
columns by strings of the form $23\ldots k$. We move through the columns of $b$  from top to bottom,
left to right. Each $+$ in $P$ (starting with the leftmost moving to the right ignoring $+$ at height $n$)
will alter $b$ as we move through the columns. Suppose the $+$ is at height $h$ in $P$.
If one encounters a spin column of type $B_n$, replace it by a column 
$12\ldots h \; h+2\ldots n\; \ol{h+1}$ (read from bottom to top). Otherwise, if one encounters a 
$\ol{1}$, replace $\ol{1}$ by $\ol{h+1}$. If one encounters a $2$, replace the string $23\ldots k$ by 
$12\ldots h\; h+2\ldots k$.

\subsection{Properties of $\pm$-diagrams for type $B_n$ and $C_n$}
\label{subsec:pm properties BC}

Let $B(\La)$ (resp. $B_{\geh_0}(\La)$) be the $C_n$ (resp. $\geh_0$)-crystal of the highest 
weight module of highest weight $\La$ and $I_0=\{1,2,\ldots,n\}$. Define
\begin{equation*}
(m_1,m_2,\ldots,m_n) = \begin{cases} 
	(2,\ldots,2,2) & \text{for type $C_n$,}\\
	(2,\ldots,2,1) & \text{for type $B_n$.}
\end{cases}
\end{equation*}
By~\cite[Theorem~5.1]{Ka:1996} with $\xi=\id$, there exists a unique injective 
map $\ol{S}:B_{\geh_0}(\La)\rightarrow B(2\La)$ such that $\ol{S}(e_ib)=e_i^{m_i}\ol{S}(b)$ 
and $\ol{S}(f_ib)=f_i^{m_i}\ol{S}(b)$ for $i\in I_0$.

\begin{lemma} \label{lem:double pm}
Let $\La=\sum_{i=1}^n k_i \La_i$, $b \in B_{\geh_0}(\La)$ be a $\{2,\ldots,n\}$-highest element,
and $P$ the corresponding $\pm$-diagram. Then the $\pm$-diagram corresponding to $\ol{S}(b)$ 
in $B(\hat{\La})$, where $\hat{\La}=\sum_{i=1}^n m_ik_i\La_i$, is obtained by doubling each column of 
$P$ together with its signs. This doubling procedure needs special treatment
when a sign is in a spin column or $0$ is at height $n$ for type $B_n$. If the former occurs, we 
replace the spin column with a full column of the same sign. If the latter occurs, we replace it with 
two columns containing $+$ and $-$.
\end{lemma}

\begin{proof}
First we prove the claim when $P$ contains only $-$ except at height $n$. At height $n$ we allow 
the sign to be $+$ or $0$. 
Let $\la_i$ be the length of the $i$-th row of the inner shape minus the
number of columns of height $n$ with $+$, and $\mu_i$ the number of $-$'s in the $i$-th row. 
If there is a spin column in type $B_n$, we regard its length or number of $-$ (if it exists) to be $1/2$. 
Let $\nu$ be the number of columns with $0$. $\nu=0$ for type $C_n$ and $0$ or $1$ for type $B_n$.
Set $M=\mu_1+\mu_2+\cdots+\mu_n$. Let $u$ be the $\{1,2,\ldots,n\}$-highest element vector of 
weight given by the outer shape of $P$. Set
\[
  \vec{a}=(1^{M+\la_1}2^{M+\la_2}\cdots(n-1)^{M+\la_{n-1}}n^{\alpha M+\nu}
  (n-1)^{\mu_1+\cdots+\mu_{n-1}}\cdots2^{\mu_1+\mu_2}1^{\mu_1}),
\]
where $\alpha=1$ for $\geh_0$ of type $C_n$ and $\alpha=2$ for $\geh_0$ of type $B_n$.
By direct calculation, one finds that $\Phi(P)=f_{\vec{a}}(u)$. Hence it is clear that $\ol{S}(\Phi(P))=
\Phi(\hat{P})$, where $\hat{P}$ is the doubled $\pm$-diagram of $P$.

Now we show the claim for the general case by induction on $N$, where $N$ is the number of $+$ 
in $P$ at height lower than $n$. If $N=0$, the statement was shown above. Suppose $N>0$. Let 
$h$ be the height of the lowest $+$ in $P$, and $P'$ the $\pm$-diagram obtained by replacing the 
leftmost $+$ at height $h$ with a box. 
Let $t$ be the rightmost column of $\Phi(P)$ such that it contains $\ol{h+1}$ or does not contain 
$h+1$. Then $\Phi(P')$ differs from $\Phi(P)$ only in the column at the position where $t$ is situated 
and there are three cases, except when $N=1$ and $P$ has a spin column with $-$. 
(This exceptional case can be checked easily and is not considered later.)
Let $t'$ be the corresponding column in $P'$. Then the three cases are as follows:
\begin{align*}
\text{(i)} \quad & t=
\begin{array}{|c|}
\hline 
k\\ \hline
\vdots\\ \hline
h+2\\ \hline
h\\ \hline
\vdots\\ \hline
1\\ \hline
\end{array},
\qquad t'=
\begin{array}{|c|}
\hline 
k\\ \hline
\vdots\\ \hline
h+2\\ \hline
h+1\\ \hline
\vdots\\ \hline
2\\ \hline
\end{array}\\[2mm]
\text{(ii)} \quad & t=
\begin{array}{|c|}
\hline 
\ol{k'}^{\hspace{-2mm}\phantom{L}}\\ \hline
k\\ \hline
\vdots\\ \hline
h+2\\ \hline
h\\ \hline
\vdots\\ \hline
1\\ \hline
\end{array},
\qquad t'=
\begin{array}{|c|}
\hline 
\ol{k'}^{\hspace{-2mm}\phantom{L}}\\ \hline
k\\ \hline
\vdots\\ \hline
h+2\\ \hline
h+1\\ \hline
\vdots\\ \hline
2\\ \hline
\end{array} \\[2mm]
\text{(iii)} \quad & t=
\begin{array}{|c|}
\hline 
\ol{h+1}^{\hspace{-2mm}\phantom{L}}\\ \hline
k\\ \hline
\vdots\\ \hline
2\\ \hline
\end{array},
\qquad t'=
\begin{array}{|c|}
\hline 
\ol{1}^{\hspace{-2mm}\phantom{L}}\\ \hline
k\\ \hline
\vdots\\ \hline
2\\ \hline
\end{array} \; .
\end{align*}
Here $k,k'\ge h+1$, and $k=h+1$ case in (i) or (ii) means that there is no letter from $h+2$ to $k$. 
To the left of $t$ we have columns of type
\[
\begin{array}{|c|}
\hline 
k\\ \hline
\vdots\\ \hline
h'+2\\ \hline
h'\\ \hline
\vdots\\ \hline
1\\ \hline
\end{array}\qquad
\text{or} \qquad
\begin{array}{|c|}
\hline 
\ol{k'}^{\hspace{-2mm}\phantom{L}}\\ \hline
k\\ \hline
\vdots\\ \hline
h'+2\\ \hline
h'\\ \hline
\vdots\\ \hline
1\\ \hline
\end{array}
\]
where $h'\ge h,k\ge h'+1,k'\ge h+1$, and to the right we have 
\[
\begin{array}{|c|}
\hline 
k\\ \hline
\vdots\\ \hline
2\\ \hline
\end{array}\qquad
\text{or} \qquad
\begin{array}{|c|}
\hline 
\ol{1'}^{\hspace{-2mm}\phantom{L}}\\ \hline
k\\ \hline
\vdots\\ \hline
2\\ \hline
\end{array} \; .
\]
Let $t^*$ be the leftmost column of $\Phi(P)$ that is to the right of $t$ and does not contain $h+1$. 
Let $n_j$ ($2\le j\le h$) be the number of boxes with letter $j$ that is weakly right of $t^*$. Setting
\begin{align*}
\vec{a}&=(12^{n_2+1}\cdots(h-1)^{n_{h-1}+1}h^{n_h+1}),\\
\vec{c}&=(h^{n_h}(h-1)^{n_{h-1}}\cdots2^{n_2}),
\end{align*}
we have $\Phi(P')=e_{\vec{c}}f_{\vec{a}}\Phi(P)$ in all three cases (i),(ii),(iii).

Next we consider the doubled $\pm$-diagrams $\hat{P}$ and $\hat{P'}$ corresponding to $P$ and 
$P'$. Then $\Phi(\hat{P'})$ differs from $\Phi(\hat{P})$ by two columns as in (i') and one column as
in (ii'), where
\begin{align*}
\text{(i')} \quad & \Phi(\hat{P}):\quad
\begin{array}{|c|c|}
\hline 
k & k\\ \hline
\vdots & \vdots\\ \hline
h+2 & h+2\\ \hline
h & h\\ \hline
\vdots & \vdots\\ \hline
1 & 1\\ \hline
\end{array},
&&\Phi(\hat{P'}):\quad
\begin{array}{|c|c|}
\hline 
k & k\\ \hline
\vdots & \vdots\\ \hline
h+2 & h+2\\ \hline
h+1 & h+1\\ \hline
\vdots & \vdots\\ \hline
2 & 2\\ \hline
\end{array}\\[2mm]
\text{(ii')} \quad & \Phi(\hat{P}):\quad 
\begin{array}{|c|}
\hline 
\ol{h+1}^{\hspace{-2mm}\phantom{L}}\\ \hline
k\\ \hline
\vdots\\ \hline
h+2\\ \hline
h\\ \hline
\vdots\\ \hline
1\\ \hline
\end{array},
&&\Phi(\hat{P'}):\quad
\begin{array}{|c|}
\hline 
\ol{1}^{\hspace{-2mm}\phantom{L}}\\ \hline
k\\ \hline
\vdots\\ \hline
h+2\\ \hline
h+1\\ \hline
\vdots\\ \hline
2\\ \hline
\end{array} \; .
\end{align*}
Let $\hat{t}^*$ be the leftmost column of $\Phi(\hat{P})$ that is to the right of the above one and 
does not contain $h+1$. The number of boxes with letter $j$ that is weakly right to $\hat{t}^*$ is given 
by $2c_j$. Calculating carefully we obtain $\Phi(\hat{P'})=e_{\vec{c'}}f_{\vec{a'}}\Phi(\hat{P})$, where 
$\vec{a'}$ and $\vec{c'}$ are sequences obtained from $\vec{a},\vec{c}$ by repeating each letter twice,
namely, if $\vec{a}=(a_1,a_2,\ldots, a_m)$, then $\vec{a'}=(a_1^2, a_2^2, \ldots, a_m^2)$. 
Therefore, we have 
\[
S(\Phi(P))=e_{\stackrel{\leftarrow}{a'}}f_{\stackrel{\leftarrow}{c'}}S(\Phi(P'))
=e_{\stackrel{\leftarrow}{a'}}f_{\stackrel{\leftarrow}{c'}}\Phi(\hat{P'})=\Phi(\hat{P}),
\]
by using the induction hypothesis.
The proof is finished.
\end{proof}

\section{Combinatorial KR crystals}
\label{sec:KR}

In this section we define the combinatorial KR crystals $V^{r,s}$. In Sections~\ref{subsec:A}
and~\ref{subsec:DBA} we review types $A_{n-1}^{(1)}$ and $D_n^{(1)}$, $B_n^{(1)}$, 
$A_{2n-1}^{(2)}$, respectively. In Section~\ref{subsec:C} we use the folding technique to construct 
$V^{r,s}$ of type $C_n^{(1)}$ as a ``virtual'' crystal inside $V^{r,s}$ of type $A_{2n+1}^{(2)}$, and 
use this in Section~\ref{subsec:A(2)D(2)} to define the combinatorial KR crystals of type 
$D_{n+1}^{(2)}$, $A_{2n}^{(2)}$ using the similarity method of Kashiwara~\cite{Ka:1996}.
The construction for exceptional nodes will be given in Section~\ref{sec:exceptional}.
In Section~\ref{subsec:phi0} we derive some properties of $\vp_0$ that will be used later.

\subsection{Combinatorial KR crystals of type $A_{n-1}^{(1)}$}
\label{subsec:A}
The existence of the KR crystal $B^{r,s}$ of type $A_{n-1}^{(1)}$
was shown in~\cite{KMN2:1992}. A combinatorial description of this crystal was 
provided by Shimozono~\cite{Sh:2002}. As a $\{1,2,\ldots, n-1\}$-crystal
\begin{equation}\label{eq:A decomp}
	B^{r,s} \cong B(s\La_r)
\end{equation}
which as a set are all semi-standard Young tableaux of rectangular shape $(s^r)$ over
the alphabet $1\prec 2\prec \cdots \prec n$. The Dynkin diagram of $A_{n-1}^{(1)}$ has a 
cyclic automorphism $i\mapsto i+1 \pmod{n}$. The action of the affine crystal operators
$f_0$ and $e_0$ is given by
\begin{equation}\label{eq:A e0}
	f_0 = \pr^{-1} \circ f_1 \circ \pr \qquad \text{and} \qquad e_0 = \pr^{-1} \circ e_1 \circ \pr,
\end{equation}
where $\pr$ is Sch\"utzenberger's promotion operator~\cite{Schuetzenberger:1972}, which
is the analogue of the cyclic Dynkin diagram automorphism on the level of crystals.
On a rectangular tableau $b \in B^{r,s}$, $\pr(b)$ is obtained from $b$ by removing all letters $n$,
adding one to each letter in the remaining tableau, using jeu-de-taquin to slide all letters up,
and finally filling the holes with $1$s.

\begin{theorem} \cite[Section 3.3]{Sh:2002}
\label{thm:A}
For type $A_{n-1}^{(1)}$, the crystal $B^{r,s}$ decomposes as in~\eqref{eq:A decomp} as a 
$\{1,2,\ldots,n-1\}$-crystal with the affine crystal action as given in~\eqref{eq:A e0}.
\end{theorem}

\subsection{Combinatorial KR crystals of type $D_n^{(1)}, B_n^{(1)}, A_{2n-1}^{(2)}$}
\label{subsec:DBA}

Let $\geh$ be of type $D_n^{(1)}$, $B_n^{(1)}$, or $A_{2n-1}^{(2)}$ with the 
underlying finite Lie algebra $\geh_0$ of type $X_n=D_n,B_n$, or $C_n$, respectively.
In this section we review the combinatorial model for KR crystals $B^{r,s}$ of type
$\geh$ as given in~\cite{OS:2008,S:2008}, where $r\neq n-1,n$ for $D_n^{(1)}$.
The case $r=n$ for $B_n^{(1)}$ treated in Lemma~\ref{lem:Vn B} in this section is new.
The cases $r=n-1,n$ for $D_n^{(1)}$ are treated in Section~\ref{subsec:exceptional D}.
The crystals for type $A_{2n-1}^{(2)}$ will be used in Section~\ref{subsec:C} in order to define 
the combinatorial KR crystals for type $C_n^{(1)}$ using the folding technique.

The Dynkin diagrams of type $D_n^{(1)},B_n^{(1)}$, and $A_{2n-1}^{(2)}$ all have
an automorphism interchanging nodes $0$ and $1$. The analogue $\sigma$ of this
automorphism on the level of crystals exists.
By construction the automorphism $\sigma$ commutes with $f_{i}$ and $e_{i}$ for $i=2,3,\ldots,n$. 
Hence it suffices to define $\sigma$ on $X_{n-1}$ highest weight elements
where $X_{n-1}$ is the subalgebra whose Dynkin diagram is obtained
from that of $X_n$ by removing node $1$. Because of the
bijection $\bij$ between $\pm$-diagrams and $X_{n-1}$-highest weight
elements as described in Section~\ref{subsec:branching}, it suffices to define the map on 
$\pm$-diagrams.

For the following description of the map $\sigD$, we assume that $1\le r\le n$ for $A_{2n-1}^{(2)}$,
$1\le r\le n-1$ for $B_n^{(1)}$ and $1\le r\le n-2$ for $D_n^{(1)}$.
Let $P$ be a $\pm$-diagram of shape $\La/\la$. Let $c_i=c_i(\la)$ be
the number of columns of height $i$ in $\la$ for all $1\le i<r$ with
$c_0=s-\la_1$. If $i\equiv r-1 \pmod{2}$, then in $P$, above each
column of $\la$ of height $i$, there must be a $+$ or a $-$.
Interchange the number of such $+$ and $-$ symbols. If $i\equiv r
\pmod{2}$, then in $P$, above each column of $\la$ of height $i$,
either there are no signs or a $\mp$ pair. Suppose there are $p_i$
$\mp$ pairs above the columns of height $i$. Change this to
$(c_i-p_i)$ $\mp$ pairs. The result is $\sigD(P)$, which has the
same inner shape $\la$ as $P$ but a possibly different outer shape.

\begin{definition} \label{def:sigma}
Let $1\le r\le n$ for $A_{2n-1}^{(2)}$, $1\le r\le n-1$ for $B_n^{(1)}$ and $1\le r\le n-2$ for $D_n^{(1)}$.
Let $b\in B^{r,s}$ and $e_{\aar} := e_{a_1} e_{a_2} \cdots e_{a_\ell}$ be such that
$e_{\aar}(b)$ is a $X_{n-1}$ highest weight crystal element. Define 
$f_{\aal}:= f_{a_\ell} f_{a_{\ell-1}} \cdots f_{a_1}$. Then 
\begin{equation} \label{eq:def sigma}
\sigma(b) := f_{\aal} \circ \bij \circ \sigD \circ \bij^{-1} \circ e_{\aar}(b).
\end{equation}
\end{definition}

As an $X_n$-crystal, $B^{r,s}$ decomposes into the following irreducible components
\begin{equation} \label{eq:classical decomp}
	B^{r,s} \cong \bigoplus_\La B(\La),
\end{equation}
with the exception of $r=n-1,n$ for $D_n^{(1)}$.
Here $B(\La)$ is the $X_n$-crystal of highest weight $\La$ and the sum runs over all dominant weights $\La$ 
that can be obtained from $s\La_r$ by the removal of vertical dominoes, where $\La_i$ are the fundamental weights 
of $X_n$.  The decomposition~\eqref{eq:classical decomp} also holds as a $\{0,2,3,\ldots,n\}$-crystal.

\begin{definition}
Let $V^{r,s}$ for $s\ge 1$ and $1\le r\le n$ for type $A_{2n-1}^{(2)}$, $1\le r\le n-1$ for type $B_n^{(1)}$, and 
$1\le r\le n-2$ for type $D_n^{(1)}$ be defined as follows. As an $X_n$ crystal
\begin{equation} \label{eq:V DBA decomp}
	V^{r,s} \cong \bigoplus_\La B(\La),
\end{equation}
where the sum runs over all dominant weights $\La$ that can be obtained from $s\La_r$ by the removal of 
vertical dominoes. The affine crystal operators $e_{0}$ and $f_{0}$ are defined as
\begin{equation} \label{eq:e0}
\begin{split}
f_{0} &= \sigma \circ f_{1} \circ \sigma,\\
e_{0} &= \sigma \circ e_{1} \circ \sigma.
\end{split}
\end{equation}
\end{definition}

Next we give a definition of $V^{n,s}$ for type $B_n^{(1)}$ using Kashiwara's similarity technique~\cite{Ka:1996}.
Note that the classical decomposition~\eqref{eq:classical decomp} is still valid in this case.

\begin{lemma}
\label{lem:Vn B}
Let $B^{n,s}$ be the $A_{2n-1}^{(2)}$-KR crystal. Then there exists a regular $B_n^{(1)}$-crystal
$V^{n,s}$ and a unique injective map $S:V^{n,s}\rightarrow B^{n,s}$ such that 
\[
S(e_ib)=e_i^{m_i}S(b),\quad S(f_ib)=f_i^{m_i}S(b)\quad\text{for }i\in I,
\]
where $(m_i)_{0\le i\le n}=(2,2,\ldots,2,1)$. Furthermore, $V^{n,s}$ decomposes 
as~\eqref{eq:V DBA decomp} as a $\{1,2,\ldots,n\}$-crystal and as a $\{0,2,\ldots,n\}$-crystal.
\end{lemma}

\begin{proof}
As a $B_n$-crystal define
\[
V^{n,s}=\bigoplus_{\bf k}B(k_\iota\La_\iota+k_{\iota+2}\La_{\iota+2}+\cdots+k_n\La_n)
\]
where the sum is over ${\bf k}\in\{(k_\iota,k_{\iota+2},\ldots,k_n)\mid
2(k_\iota+k_{\iota+2}+\cdots+k_{n-2})+k_n=s\}$ with $\iota\equiv n$\,(mod $2$),\;$\iota=0\text{ or }1$,
and $B(\La)$ is the highest weight $B_n$-crystal of highest weight $\La$, where we identify $\La_0=0$.
Then by~\cite[Theorem 5.1]{Ka:1996} with $\xi=\id$ and $(m_i)_{i\in\{1,2,\ldots,n\}}$, one realizes that there exists 
a unique injective map
\[
\ol{S}:V^{n,s}\longrightarrow\bigoplus_{\bf k}
B_{C_n}(2k_\iota\La_\iota+2k_{\iota+2}\La_{\iota+2}+\cdots+2k_{n-2}\La_{n-2}+k_n\La_n)
\]
such that $\ol{S}(e_ib)=e_i^{m_i}\ol{S}(b)$ and $\ol{S}(f_ib)=f_i^{m_i}\ol{S}(b)$ for $i\in\{1,2,\ldots,n\}$.
Note that the RHS is contained in $B^{n,s}$. Since $\sigma$ is defined on $B^{n,s}$, to finish the proof 
it suffices to show that $\ol{S}(V^{n,s})$ is closed under $\sigma$. This is clear from 
Lemma~\ref{lem:double pm}.

By construction, $V^{n,s}$ decomposes as~\eqref{eq:V DBA decomp} as a $\{1,2,\ldots,n\}$-crystal.
Since $e_0$ is defined using $\sigma$, the same decomposition holds as a 
$\{0,2,\ldots,n\}$-crystal.
\end{proof}

Let us now state \cite[Lemma 5.1]{S:2008} since we will need it several times in the 
sequel. The lemma describes the action of $e_1$ on $\{3,4,\ldots,n\}$ highest weight
elements which are in one-to-one correspondence with pairs of $\pm$-diagrams $(P,p)$,
where the inner shape of $P$ is the outer shape of $p$.

The operator $e_1$ either changes a 2 into a 1 or a $\ab$ into a $\bb$
in $b$ corresponding to the pair of $\pm$-diagrams $\Psi(P,p)$. On the level of $(P,p)$ this means 
that either a $+$ from $p$ transfers to $P$,
or a $-$ moves from $P$ to $p$. To describe the precise action of $e_1$ on $(P,p)$ perform
the following algorithm:
\begin{enumerate}
\item Successively run through all $+$ in $p$ from left to right and, if possible, pair it with 
the leftmost yet unpaired $+$ in $P$ weakly to the left of it.
\item Successively run through all $-$ in $p$ from left to right and, if possible, pair it with
the rightmost  yet unpaired $-$ in $P$ weakly to the left.
\item Successively run through all yet unpaired $+$ in $p$ from left to right and, if possible,
pair it with the leftmost yet unpaired $-$ in $p$.
\end{enumerate}

\begin{lemma} \cite[Lemma 5.1]{S:2008} \label{lem:e1 action}
If there is an unpaired $+$ in $p$,  $e_1$ moves the rightmost unpaired $+$ in $p$ to $P$. 
Else, if there is an unpaired $-$ in $P$, $e_1$ moves the leftmost unpaired $-$ in $P$ to $p$.
Else $e_1$ annihilates $(P,p)$.
\end{lemma}

\subsection{Combinatorial KR crystals for $C_n^{(1)}$}
\label{subsec:C}

In this section we define a combinatorial crystal $V^{r,s}$ of type $C_n^{(1)}$.
For $r\neq n$, this combinatorial crystal is isomorphic to the KR crystal $B^{r,s}$ (see 
Theorem~\ref{thm:uniqueness}). The KR crystal for $r=n$ is treated in 
Section~\ref{subsec:exceptional CD}.
We realize $V^{r,s}$ as a ``virtual'' crystal inside the ambient crystal 
$\Vh^{r,s}=B^{r,s}_{A_{2n+1}^{(2)}}$ of type $A_{2n+1}^{(2)}$. Let $I=\{0,1,\ldots,n\}$
be the index set for the Dynkin diagram of type $C_n^{(1)}$ and $\Ih=\{0,1,\ldots,n+1\}$ be the index
set of the Dynkin diagram of type $A_{2n+1}^{(2)}$.  We denote the
crystal operators of $V^{r,s}$ by $e_i$ and $f_i$, and the crystal operators of the ambient 
crystal $\Vh^{r,s}$ by $\eh_i$ and $\fh_i$. In this case
\begin{equation} \label{eq:operators}
	e_i = \begin{cases} \eh_0\eh_1 & \text{if $i=0$}\\
	                                    \eh_{i+1} & \text{if $1\le i\le n$}
	           \end{cases}
\quad \text{and} \quad
	 f_i = \begin{cases} \fh_0\fh_1 & \text{if $i=0$}\\
	                                    \fh_{i+1} & \text{if $1\le i\le n$.}
	           \end{cases}
\end{equation}
The analogue of the Dynkin diagram automorphism of $A_{2n+1}^{(2)}$, which interchanges nodes
$0$ and $1$, on the level of the crystal $B^{r,s}_{A_{2n+1}^{(2)}}$ is denoted by $\sigma$.

\begin{definition}
Define $V^{r,s}$ to be the subset of elements $b\in\Vh^{r,s}$ that are invariant under $\sigma$, 
namely $\sigma(b)=b$, together with the operators $e_i$ and $f_i$ of~\eqref{eq:operators}.
\end{definition}

\begin{lemma} \label{lem:closed}
$V^{r,s}$ is closed under the operators $e_i$ and $f_i$ for $i\in I$.
\end{lemma}

\begin{proof}
For $i\in \Ih\setminus \{0,1\}$, the operators $\eh_i$ and $\fh_i$ commute with
$\sigma$. Hence for $b\in V^{r,s}$ such that $e_i(b)\neq \emptyset$ we have
\begin{equation*}
	\sigma \circ e_i(b) = e_i\circ \sigma(b) = e_i(b) \quad \text{for $i\in I\setminus \{0\}$}
\end{equation*}
and similarly for $f_i$. Note also that $\eh_0=\sigma \circ \eh_1 \circ \sigma$ and 
$\fh_0=\sigma \circ \fh_1 \circ \sigma$. Hence for an element $b\in V^{r,s}$ such that
$e_0(b)\neq \emptyset$ we have
\begin{multline*}
	\sigma\circ e_0(b) = \sigma\circ \eh_0\circ \eh_1(b) = \sigma\circ \eh_1\circ \eh_0(b)
	= \sigma\circ \eh_1\circ \sigma \circ \eh_1\circ \sigma(b)\\
	=  (\sigma\circ \eh_1\circ \sigma) \circ \eh_1(b) = \eh_0 \circ \eh_1(b) = e_0(b)
\end{multline*}
and similarly for $f_0$.
\end{proof}

\begin{lemma} \label{lem:decomp}
As a $\{1,2,\ldots,n\}$-crystal
$$
	V^{r,s} \cong \bigoplus_\Lambda B(\Lambda)
$$
where the sum is over all $\Lambda$ obtained from $s\Lambda_r$ by removing horizontal
dominoes and $B(\Lambda)$ is a highest weight $C_n$-crystal of highest weight $\Lambda$.
\end{lemma}

\begin{proof}
As a $\{1,2,\ldots,n\}$-crystal $V^{r,s}$ is isomorphic to the connected components of $\sigma$ 
invariant elements of the $\{2,3,\ldots,n+1\}$-crystal $\Vh^{r,s}$ of type $A_{2n+1}^{(2)}$. Since 
$\eh_i$ and $\fh_i$ commute with $\sigma$ for $i\in \Ih\setminus \{0,1\}$, the crystal of $\sigma$ 
invariant elements of the $\{2,3,\ldots,n+1\}$-crystal $\Vh^{r,s}$ of type $A_{2n+1}^{(2)}$ is a disjoint 
union of $C_n$-crystals. The $\{2,3,\ldots,n+1\}$-highest weight elements in $\Vh^{r,s}$ are in
one-to-one correspondence with $\pm$-diagrams.  Hence to find the decomposition we need
to list all $\sigma$-invariant $\pm$-diagrams. Let $P$ be a $\pm$-diagram for $\Vh^{r,s}$ of
type $A_{2n+1}^{(2)}$. Let $\lambda$ be the outer shape of $P$ and $\Lambda$ be the inner
shape of $P$. For $P$ to be $\sigma$-invariant, for any row $i$ of $\Lambda$ with $r-i\equiv 1
\pmod{2}$, the cells above the cells in row $i$ of $\Lambda$ must contain the same number of
$+$ and $-$ signs in $P$. Similarly, for rows $i<r$ in $\Lambda$ with $r-i\equiv 0 \pmod{2}$,
the cells above the cells in row $i$ of $\Lambda$ must contain the same number of
$\pm$ pairs as no additional cells in $P$. This shows in particular that the rows of $\Lambda$
have the same parity as $s$, and hence $\Lambda$ can be obtained from $s\Lambda_r$
by removal of horizontal dominoes. As a $C_n$-weight, the weight of the $\{2,3,\ldots,n+1\}$-
highest weight element corresponding to $P$ is precisely $\Lambda$. This proves the claim.
\end{proof}

\begin{lemma} \label{lem:decomp1}
As a $\{0,1,\ldots,n-1\}$-crystal
$$
	V^{r,s} \cong \bigoplus_\Lambda B(\Lambda)
$$
where the sum is over all $\Lambda$ obtained from $s\Lambda_r$ by removing horizontal
dominoes and $B(\Lambda)$ is a highest weight $C_n$-crystal of highest weight $\Lambda$.
\end{lemma}

\begin{proof}
Proposition 3.2.1 of~\cite{NS:2006} applies to the situation, where 
$\mathfrak{g}=A_{2n+1}^{(2)}$, $\hat{\mathfrak{g}}=C_n^{(1)}$, $\omega$ is
the Dynkin diagram automorphism interchanging nodes $0$ and $1$ (this Proposition
applies even though $\omega(0)\neq 0$ in this case), and $B=B^{r,s}$ of type $A_{2n+1}^{(2)}$.
Then the Proposition ensures that the fixed point subset $B^\omega$ (or $V^{r,s}$ in our notation)
with crystal operators as in~\eqref{eq:operators} is a regular $U_q'(\hat{\mathfrak{g}})$-crystal.

The Weyl group of type $C_n$ contains an element that maps
$\Lambda_j-\Lambda_0$ to $\Lambda_{n-j}-\Lambda_n$ for all $1\le j\le n$.
Under this map, the highest weight elements in $V^{r,s}$ as a $\{1,2,\ldots,n\}$-crystal
map to highest weight elements in $V^{r,s}$ as a $\{0,1,\ldots,n-1\}$-crystal.
Hence the decomposition of the claim follows from Lemma~\ref{lem:decomp}.
\end{proof}

\begin{lemma} \label{lem:double}\mbox{}
\begin{enumerate}
\item
Let $b\in V^{r,s}$ be a $\{2,\ldots,n\}$-highest element corresponding to the $\pm$-diagram $P$.
Suppose the number of $\mp,+,-$ or $\cdot$ on the columns of the inner shape of the same height
smaller than $r$ is always even. Then the $\pm$-diagram corresponding to $e_0(b)$ also has the 
same property.
\item
Let $\La, \La'$ be weights appearing in the decomposition of $V^{r,s}$ as in Lemma~\ref{lem:decomp}
such that $\La'$ is obtained from $\La$ by adding a horizontal domino. Then there exists a
$\{2,\ldots,n\}$-highest weight element $b\in B(\La)$ such that $e_0(b)\in B(\La')$.
\end{enumerate}
\end{lemma}

\begin{proof}
Let us first prove (1).
Inside the ambient crystal $\Vh^{r,s}$, the element $b$ is a $\{3,\ldots,n+1\}$
highest weight vector. It corresponds to a pair of $\pm$-diagrams $(P',P)$, where $P$ is
the same as in the statement of the lemma and $P'$ is the $\pm$-diagram corresponding
to the highest weight vector in the component of $b$ obtained as described in the proof of 
Lemma~\ref{lem:decomp}. The crystal operator $e_0$ in $V^{r,s}$ corresponds to
$\eh_1 \eh_0= \eh_1 \sigma \eh_1 \sigma = \eh_1 \sigma \eh_1$ in $\Vh^{r,s}$, where in
the last step we used that $b\in V^{r,s}$ is invariant under $\sigma$.
By Lemma~\ref{lem:e1 action}, $\eh_1$ either moves a $+$ from $P$ to $P'$, or a $-$ from
$P'$ to $P$. In this case $\sigma$ on $\eh_1(P')$ either has one extra $-$, or one fewer $+$,
respectively. Since by Lemma~\ref{lem:closed}, $V^{r,s}$ is closed under $e_0$ and the number
of $+$ and $-$ is balanced in $\pm$-diagrams corresponding to $\{2,\ldots,n\}$-highest weight
vectors, $e_0(P')$ must have this property. This implies that $e_0(P)$ either has two less $+$
in the same row or two extra $-$ in the same row (since the shape can
only change by horizontal dominoes). This proves the claim.

For the proof of (2), let $P'$ be as before with inner shape $\La$. Suppose $\La'$ has two more 
boxes than $\La$ in row $j$. Let $P$ be the $\pm$-diagram with outer shape $\La$ corresponding 
to the partition $\la=(\la_1,\ldots,\la_n)$ such that row $n$ is filled with $-$ and each row $\la_i$ for 
$n>i\neq j-1$ is filled with $(\la_i-\la_{i+1})/2$ minus signs. Then by the previous description of $e_0$ 
using Lemma~\ref{lem:e1 action}, the $\pm$-diagram corresponding to $e_0(b)$ has outer shape 
$\La'$ and two more minus signs in row $j$.
\end{proof}

\subsection{Combinatorial KR crystals for $A_{2n}^{(2)},D_{n+1}^{(2)}$} 
\label{subsec:A(2)D(2)}
We use the similarity technique of Kashiwara developed in \cite{Ka:1996} to define combinatorial
models of KR crystals for $A_{2n}^{(2)},D_{n+1}^{(2)}$ for nonexceptional nodes from those for 
$C_n^{(1)}$. In this subsection $V^{r,s}$ stands for the $C_n^{(1)}$-crystal defined in the previous subsection.
For types $A_{2n}^{(2)},D_{n+1}^{(2)}$ define positive integers $m_i$ for $i\in I$
as follows:
\begin{equation}
\label{eq:m}
(m_0,m_1,\ldots,m_{n-1},m_n) = \begin{cases}
       (1,2,\ldots,2,2) & \text{ for }A_{2n}^{(2)},\\
       (1,2,\ldots,2,1) & \text{ for }D_{n+1}^{(2)}.
\end{cases}
\end{equation}

The next theorem ensures the existence of crystals $V_\geh^{r,s}$ for type $A_{2n}^{(2)}$ and
$D_{n+1}^{(2)}$. They are constructed from type $C_n^{(1)}$ crystals by an affine extension of 
$\ol{S}$ of Section~\ref{subsec:pm properties BC}.

\begin{theorem}
\label{th:embedding}
For $1\le r\le n$ for $\geh=A_{2n}^{(2)}$, $1\le r<n$ for $\geh=D_{n+1}^{(2)}$ and $s\ge1$,
there exists a $\geh$-crystal $V^{r,s}_\geh$ and a unique injective map
$S:V^{r,s}_\geh\longrightarrow V^{r,2s}$ such that
\[
  S(e_ib)=e_i^{m_i}S(b),\quad S(f_ib)=f_i^{m_i}S(b) \qquad \text{for $i \in I$.}
\]
\end{theorem}

\begin{proof}
Let $\La$ be a dominant integral weight of $\geh_0$ where $\geh_0$ is $C_n,B_n$ for 
$\geh=A_{2n}^{(2)},D_{n+1}^{(2)}$, respectively. We assume if $\geh=D_{n+1}^{(2)}$, the coefficient 
of the $n$-th fundamental weight of $\La$ is zero. As before let $B(\La)$ (resp. $B_{\geh_0}(\La)$) 
be the $C_n$ (resp. $\geh_0$)-crystal of the highest weight module of highest weight $\La$. 
We saw that by~\cite[Theorem~5.1]{Ka:1996}, there exists a unique injective map 
$\ol{S}:B_{\geh_0}(\La)\rightarrow B(2\La)$ such that $\ol{S}(e_ib)=e_i^{m_i}\ol{S}(b)$ and 
$\ol{S}(f_ib)=f_i^{m_i}\ol{S}(b)$ for $i\in I_0=\{1,2,\ldots,n\}$.

Now define $V_{\geh}^{r,s}$, as a $\geh_0$-crystal, by
\begin{equation} \label{eq:V decomp}
   V_{\geh}^{r,s}=\bigoplus_{\La}B_{\geh_0}(\La),
\end{equation}
where the sum is over all $\La$ obtained from $s\La_r$ by removing nodes. Then from the above
we have an injective map
$S:V_{\geh}^{r,s}\rightarrow V^{r,2s}$. To make $V_{\geh}^{r,s}$ a $\geh$-crystal, it
suffices to introduce the 0-action on $V_{\geh}^{r,s}$. Take an element $b\in V_{\geh}^{r,s}$. One
can assume $b$ is $\{2,\ldots,n\}$-highest, hence $b$ corresponds to a $\pm$-diagram $P$. Then
$S(b)$ corresponds to a $\pm$-diagram $\tilde{P}$, where the number of $\mp,+,-$ or $\cdot$
on the columns of the inner shape of the same height is doubled from that of $P$ by 
Lemma~\ref{lem:double pm}. By Lemma~\ref{lem:double}, $e_0S(b)$ also has the same property, 
so there exists a $b'\in V_{\geh}^{r,s}$ such that $S(b')=e_0S(b)$. Then one can define $e_0b=b'$ on
$V_{\geh}^{r,s}$ and we have $S(e_0b)=e_0S(b)$ by definition. The case $f_0$ is similar.
\end{proof}

\begin{lemma}
\label{lem:decomp2}
For $1\le r\le n$ for $\geh=A_{2n}^{(2)}$, $1\le r<n$ for $\geh=D_{n+1}^{(2)}$ and $s\ge1$,
\begin{itemize}
\item[(1)] as a $\{1,2,\ldots,n\}$-crystal
\[
V_{\geh}^{r,s}\simeq\bigoplus_{\La}B_{\geh_0}(\La)
\]
where the sum is over all $\La$ obtained from an $r\times s$ rectangle by removing boxes and 
$B_{\geh_0}(\La)$ is a highest weight $\geh_0$-crystal of highest weight $\La$;
\item[(2)] as a $\{0,1,\ldots,n-1\}$-crystal
\[
V_{\geh}^{r,s}\simeq\bigoplus_{\La}B_{B_n}(\La)
\]
where the sum is over all $\La$ obtained from an $r\times s$ rectangle by removing horizontal dominoes
for $\geh=A_{2n}^{(2)}$ and boxes for $\geh=D_{n+1}^{(2)}$, and $B_{B_n}(\La)$
is a highest weight $B_n$-crystal of highest weight $\La$.
\end{itemize}
\end{lemma}

\begin{proof}
Decomposition (1) is true by construction, see the proof of Theorem \ref{th:embedding} and
in particular~\eqref{eq:V decomp}. Once (1) is established, the decomposition
as a $\{0,1,\ldots,n-1\}$-crystal is unique. See e.g. Section 6.2 of~\cite{HKOTT:2002} for this result.
\end{proof}

\subsection{Some properties of $\vp_0$}
\label{subsec:phi0}
In this section $V^{r,s}$ is the crystal of type $C_n^{(1)}$, $A_{2n}^{(2)}$ or $D_{n+1}^{(2)}$
as defined in Sections~\ref{subsec:C} and \ref{subsec:A(2)D(2)}.

\begin{proposition} \label{prop:phi0}
Let $b\in V^{r,s}$ be a $\{2,\ldots,n\}$-highest weight element and $P$ the corresponding
$\pm$-diagram.
\begin{enumerate}
\item
Suppose that $P$ contains only empty columns $\cdot$ or $\ve$ ($\ve=+$ or $-$).  Let $\la$ be the 
inner shape of $P$. We move through the columns of $\la$ with height less than $r$ from left to right.
(We consider there are $s-\mbox{width}(b)$ columns of height $0$.) Delete all $(\cdot,\ve)$
pairs successively until we have some $\ve$'s followed by some $\cdot$'s. Then, setting
$m_0=2$ for $C_n^{(1)}$ and $m_0=1$ for $A_{2n}^{(2)}$ and $D_{n+1}^{(2)}$:
\begin{enumerate}
\item
If $\ve=+$, $\vp_0(b)$ is the number of the remaining $\cdot$ divided by $m_0$. 
\item
If $\ve=-$, $\vp_0(b)$ is the number of the remaining $-$ divided by $m_0$ plus $(s-c_r)/m_0$, 
where $c_r$ is the number of columns in $\la$ of height $r$.
\end{enumerate}
\item Suppose that the outer shape of $P$ has at least one column of height strictly less than $n-1$,
but $P$ does not contain any $+$ in columns of height strictly less than $n-1$. Then $\vp_0(b)>0$.
\end{enumerate}
\end{proposition}

\begin{proof}
We first treat case (1) for type $C_n^{(1)}$.
As in the proof of Lemma~\ref{lem:double}, inside the ambient crystal $\Vh^{r,s}$, the element $b$
is a $\{3,\ldots,n+1\}$-highest weight vector and corresponds to a pair of $\pm$-diagrams $(P',P)$,
where $P$ is the same as in the statement of the lemma and $P'$ is the $\pm$-diagram corresponding
to the highest weight vector in the component of $b$ which is invariant under $\sigma$.
Since $f_0$ corresponds to $\fh_0 \fh_1$ in $\Vh^{r,s}$, and $\fh_0$ and $\fh_1$ commute,
it suffices to determine $\widehat{\varphi}_1$ of the element corresponding to the pair of
$\pm$-diagrams $(P',P)$ to obtain $\varphi_0(b)$. We again employ Lemma~\ref{lem:e1 action}.

First suppose that $\ve=+$. Recall that this means that $P$ does not contain any $-$. After bracketing
the $+$ in $P$ with the $+$ in $P'$, $\varphi_0$ is determined by the number of unbracketed $+$
in $P'$ by Lemma~\ref{lem:e1 action}. Recall that if $b$ sits in the $C_n$ component of highest
weight $\Lambda$, then $P'$ is obtained by adding above each column in $\Lambda$ either the
same number of $+$ and $-$, or the same number of $\mp$ and $\cdot$. With this, it is not hard to
check that the rule stated in the lemma coincides with the bracketing arguments of
Lemma~\ref{lem:e1 action}.

Next suppose that $\ve=-$. In this case $P$ only contains $-$ and after bracketing the $-$ in $P$
with the $-$ in $P'$, $\vp_0$ is obtained by the number of $+$ in $P'$ (which equals
$(s-c_r)/2$) plus the number of unbracketed $-$ in $P$ (which by similar arguments as in the
case $\ve=+$ is equal to half of the number of unbracked $-$ by the rule of the lemma).

The proof for type $A_{2n}^{(2)}$ and $D_{n+1}^{(2)}$ is now immediate by 
Lemma~\ref{lem:double pm}.

For case (2), we consider again the pair $(P',P)$ as before. Since there is at least one column
of height strictly less than $n-1$ and there are no $+$s in $P$ in columns of height strictly less
than $n-1$, we conclude by Lemma~\ref{lem:e1 action} that there is at least one unbracketed 
$+$ in $P'$ and hence $\vp_0(b)>0$.
\end{proof}

\section{Uniqueness theorem}
\label{sec:uniqueness}
Recall that $B^{r,s}$ denotes the KR crystal and $V^{r,s}$ denotes the combinatorial crystals
defined in Section~\ref{sec:KR}.
In this section we show that given certain decompositions of $B^{r,s}$ into subcrystals,
$B^{r,s}$ is uniquely determined, thereby proving that $V^{r,s}\cong B^{r,s}$.
In Section~\ref{subsec:unique DBA} we review the proof of~\cite{OS:2008} for types $D_n^{(1)}$, 
$B_n^{(1)}$, $A_{2n-1}^{(2)}$ and extend it to $r=n$ of type $B_n^{(1)}$.
For type $C_n^{(1)}$, $D_{n+1}^{(2)}$, $A_{2n}^{(2)}$ the uniqueness proof is given in
Theorem~\ref{thm:uniqueness} in Section~\ref{subsec:unique CDA}.

\subsection{Uniqueness for types $D_n^{(1)}$, $B_n^{(1)}$, $A_{2n-1}^{(2)}$}
\label{subsec:unique DBA}
Let $V^{r,s}$ be the combinatorial crystal of type $D_n^{(1)}$, $B_n^{(1)}$, or $A_{2n-1}^{(2)}$
as defined in Section~\ref{subsec:DBA}.
By definition of $V^{r,s}$ we know that for $1\le r\le n-2$ for type $D_n^{(1)}$ and for
$1\le r\le n$ for types $B_n^{(1)}$ and $A_{2n-1}^{(2)}$ there exist the isomorphisms
\begin{equation*}
\begin{split}
	\Psi_0 : V^{r,s} \simeq B^{r,s} & \qquad \text{as an isomorphism of $\{1,2,\ldots,n\}$-crystals,}\\
	\Psi_1 : V^{r,s} \simeq B^{r,s} & \qquad \text{as an isomorphism of $\{0,2,\ldots,n\}$-crystals.}
\end{split}
\end{equation*}

\begin{theorem} \label{thm:DBA}
Let $s\ge 1$ and $1\le r\le n$ for type $A_{2n-1}^{(2)}$, $1\le r\le n-1$ for type $B_n^{(1)}$, and 
$1\le r\le n-2$ for type $D_n^{(1)}$. Then $\Psi_0(b)=\Psi_1(b)$ for all $b\in V^{r,s}$, and hence
there exists a unique $I$-crystal isomorphism $\Psi:V^{r,s} \cong B^{r,s}$.
\end{theorem}

We review the proof of~\cite[Proposition 6.1]{OS:2008} of this theorem here, with special attention to the
the case $r=n$ for type $A_{2n-1}^{(2)}$. This result is used in Theorem~\ref{thm:B} below for the uniqueness 
proof for $r=n$ of type $B_n^{(1)}$. Let us first recall a Lemma and Remark that is needed for the proof of
Theorem~\ref{thm:DBA}.

\begin{lemma} \label{lem:containment}
Let $b\in V^{r,s}$ be an $X_{n-2}$-highest weight vector corresponding to the tuple of $\pm$-diagrams $(P,p)$ 
where all columns of height less than $n-1$ in $p$ are empty and all columns of height $n-1$ in $p$ contain $-$.
Assume that $\ve_0(b),\ve_1(b)>0$. Then $\is(b)$ is strictly contained in 
$\is(e_{0}(b))$, $\is(e_{1}(b))$, and $\is(e_{0}e_{1}(b))$.
\end{lemma}

\begin{proof}
By assumption $p$ does not contain any $+$ and $e_{1}$ is defined. Hence by Lemma~\ref{lem:e1 action},
$e_{1}$ moves a $-$ in $P$ to $p$. This implies that the inner shape of $b$ is strictly
contained in the inner shape of $e_{1}(b)$.

The involution $\sigma$ does not change the inner shape of $b$ (only the outer shape).
By the same arguments as before, the inner shape of $b$ is strictly contained in
the inner shape of $e_{1} \sigma(b)$. Since $\sigma$ does not change the inner shape, this 
is still true for $e_{0}(b) = \sigma e_{1} \sigma(b)$.

Now let us consider $e_{0} e_{1}(b)$. For the change in inner shape
we only need to consider $e_{1} \sigma e_{1}(b)$. By the same arguments as before,  $e_{1}$
moves a $-$ from $P$ to $p$ and $\sigma$ does not change the inner shape. The next
$e_{1}$ will move another $-$ in $\sigma e_{1}(b)$ to $p$. Hence $p$ will
have grown by two $-$, so that the inner shape of $e_{1} \sigma e_{1}(b)$
is increased by two boxes.
\end{proof}

\begin{remark} \label{rem:psi weight}
Note that $\Psi_0$ and $\Psi_1$ preserve weights, that is,
$\wt(b)=\wt(\Psi_0(b))=\wt(\Psi_1(b))$ for all $b\in V^{r,s}$.
This is due to the fact that if all but one coefficient $m_j$ are known for a weight
$\Lambda=\sum_{j=0}^n m_j \La_j$, then the missing $m_j$ is also determined by the
condition that the weights in KR crystals are of level zero.
\end{remark}

\begin{proof}[Proof of Theorem~\ref{thm:DBA}]
If $\Psi_0(b)=\Psi_1(b)$ for a $b$ in a given $X_{n-1}$-component $\mathcal{C}$, then 
$\Psi_0(b')=\Psi_1(b')$ for all $b'\in \mathcal{C}$ since $e_{i}\Psi_0(b')=\Psi_0(e_{i}b')$ and 
$e_{i}\Psi_1(b')=\Psi_1(e_{i}b')$ for $i\in J=\{2,3,\ldots,n\}$. Hence it suffices to prove 
$\Psi_0(b)=\Psi_1(b)$ for only one element $b$ in each $X_{n-1}$-component $\mathcal{C}$. 
We are going to establish the theorem for $b$ corresponding to the pairs of $\pm$-diagrams $(P,p)$ 
where all columns of $p$ of height smaller than $n-1$ are empty and the columns of height $n-1$ 
are filled with $-$ ($p$ can only contain columns of height $n-1$ for 
type $A_{2n-1}^{(2)}$ and $r=n$). Note that this is an $X_{n-2}$-highest
weight vector, but not necessarily an $X_{n-1}$-highest weight vector.

We proceed by induction on $\is(b)$ by containment. First suppose that both $\ve_0(b),\ve_1(b)>0$.
By Lemma~\ref{lem:containment}, the inner shape of $e_{0} e_{1}b$, $e_{0}b$, and $e_{1}b$ is bigger
than the inner shape of $b$, so that by induction hypothesis
$\Psi_0(e_{0} e_{1}b)=\Psi_1(e_{0} e_{1}b)$, $\Psi_0(e_{0}b)=\Psi_1(e_{0}b)$, and 
$\Psi_0(e_{1}b)=\Psi_1(e_{1}b)$. Therefore we obtain
\begin{multline*}
	e_{0} e_{1} \Psi_0(b) = e_{0} \Psi_0(e_{1} b) = e_{0} \Psi_1(e_{1} b) = \Psi_1(e_{0} e_{1} b)
	= \Psi_0(e_{0} e_{1} b)\\
	 = e_{1} \Psi_0(e_{0} b) = e_{1} \Psi_1(e_{0} b) = e_{1} e_{0} \Psi_1(b).
\end{multline*}
This implies that $\Psi_0(b)=\Psi_1(b)$.

Next we need to consider the cases when $\ve_0(b)=0$ or $\ve_1(b)=0$, which comprises the base
case of the induction. Let us first treat the case $\ve_1(b)=0$.
Recall that all columns of $p$ of height smaller than $n-1$ are empty and columns of height $n-1$ are
filled with $-$. Hence it follows from the description of the action of $e_{1}$ of Lemma~\ref{lem:e1 action},
that $\ve_1(b)=0$ if and only if all columns of height smaller than $n$ in $P$ are either empty or contain $+$ and
the columns of height $n$ contain either $+$ or $-$.
\begin{quote}
\textbf{Claim.} $\Psi_0(b)=\Psi_1(b)$ for all $b$ corresponding to the pair of $\pm$-diagrams
$(P,p)$ where all columns of height smaller than $n$ in $P$ are either empty or contain $+$, and
all columns of $p$ of height smaller than $n-1$ are empty and columns of height $n-1$ are filled with $-$.
\end{quote}
The claim is proved by induction on $k+k_-$, where $k$ is the number of empty columns
in $P$ of height strictly smaller than $r$ and $k_-$ is the number of $-$ in $P$. For $k+k_-=0$ the claim is true 
by weight considerations.
Now assume the claim is true for all $0\le k'<k+k_-$ and we will establish the claim for $k+k_-$.
Suppose that $\Psi_1(b)=\Psi_0(\tilde{b})$ where $\tilde{b}\neq b$. By weight considerations 
$\tilde{b}$ must correspond to a pair of $\pm$-diagrams $(\tilde{P},p)$, where $\tilde{P}$ has the same number
of columns containing only $+$ or $-$ as $P$ and some of the empty columns of $P$ of
height $h$ strictly smaller than $r$ could be replaced by columns of height $h+2$ containing
$\mp$. Denote by $k_+$ the number of columns of height strictly less than $n$ in $P$ containing $+$. Then
\begin{equation*}
	m := \ve_0(b) = k_+ + k,
\end{equation*}
since under $\sigma$ all empty columns in $P$ become columns with $\mp$ and columns
containing $+$ become columns with $-$. By Lemma~\ref{lem:e1 action}, then $e_{1}$ acts on
$(\sigD(P),p)$ as often as there are minus signs in $\sigD(P)$ of height less than $n$, which is $k_++k$.
Set $\hat{b}=e_{1}^a\tilde{b}$, where $a>0$ is the number of columns in $\tilde{P}$ containing
$\mp$ plus the number of columns of $\tilde{P}$ of height strictly less than $n$ containing $-$. If $(\hat{P},\hat{p})$
denotes the tuple of $\pm$-diagrams associated to $\hat{b}$, then compared to $(\tilde{P},p)$ all $-$ from the 
$\mp$ pairs in $\tilde{P}$ and all $-$ in columns of height less than $n$ moved to $p$.
Note that the induction variable for $\hat{P}$ is $k+k_- -a<k+k_-$, so that by induction
hypothesis $\Psi_0(\hat{b})=\Psi_1(\hat{b})$. Hence
\begin{equation} \label{eq:eq string}
	\Psi_1(b) = \Psi_0(\tilde{b}) = \Psi_0(f_{1}^a \hat{b}) = f_{1}^a \Psi_0(\hat{b}) 
	= f_{1}^a \Psi_1(\hat{b}).
\end{equation}
Note that 
\begin{equation*}
	\ve_0(\hat{b}) = \ve_0(\tilde{b}) = m-a<m.
\end{equation*}
Hence
\begin{equation*}
\begin{split}
	&e_{0}^m \Psi_1(b) = \Psi_1(e_{0}^m b) \neq \emptyset\\
	\text{but} \qquad 
	&e_{0}^m f_{1}^a \Psi_1(\hat{b}) = f_{1}^a \Psi_1(e_{0}^m \hat{b}) = \emptyset
\end{split}
\end{equation*}
which contradicts~\eqref{eq:eq string}. This implies that we must have $\tilde{b}=b$ proving the claim.

The case $\varepsilon_0(b)=0$ can be proven in a similar fashion to the case $\varepsilon_1(b)=0$. 
Using the explicit action of $\sigD$ on $P$ and Lemma~\ref{lem:e1 action}, it follows that $\varepsilon_0(b)=0$
if and only if all columns of $P$ of height strictly less than $n$ contain either $-$ or $\mp$ pairs.
\begin{quote}
\textbf{Claim.} $\Psi_0(b)=\Psi_1(b)$ for all $b$ corresponding to the pair of $\pm$-diagrams
$(P,p)$ where all columns of $P$ of height strictly less than $n$ contain either $-$ or $\mp$ pairs, and
all columns of $p$ of height smaller than $n-1$ are empty and columns of height $n-1$ are filled with $-$.
\end{quote}
By induction on the number of $\mp$ pairs plus the number of $+$ at height $n$ in $P$, this claim can be 
proven similarly as before (using the fact that $\sigD$ changes columns with $-$ into columns with $+$ and 
columns with $\mp$ pairs into empty columns).
\end{proof}

\begin{theorem}
\label{thm:B}
For type $B_n^{(1)}$ with $V^{n,s}$ as in Lemma~\ref{lem:Vn B}, we have $V^{n,s}\cong B^{n,s}$.
\end{theorem}
\begin{proof}
The proof follows in the same way as the proof of Theorem~\ref{thm:DBA} for type $A_{2n-1}^{(2)}$, 
keeping in mind that $\pm$-diagrams in $V^{n,s}$ are characterized using 
Lemma~\ref{lem:double pm}.
\end{proof}

\subsection{Uniqueness for types $C_n^{(1)}$, $D_{n+1}^{(2)}$, $A_{2n}^{(2)}$}
\label{subsec:unique CDA}
In this section $V^{r,s}$ is the combinatorial crystal of type $C_n^{(1)}$ of Section~\ref{subsec:C}
or of type $A_{2n}^{(2)}$ or $D_{n+1}^{(2)}$ as defined in Section~\ref{subsec:A(2)D(2)}.
By Lemmas~\ref{lem:decomp}, \ref{lem:decomp1} and~\ref{lem:decomp2}, $V^{r,s}$ is isomorphic 
to $B^{r,s}$ as an $\{1,2,\ldots,n\}$-crystal and as an $\{0,1,\ldots,n-1\}$-crystal. We define the 
isomorphisms
\begin{equation*}
\begin{split}
	\Psi_0 : V^{r,s} \simeq B^{r,s} & \qquad \text{as an isomorphism of $\{1,2,\ldots,n\}$-crystals,}\\
	\Psi_n : V^{r,s} \simeq B^{r,s} & \qquad \text{as an isomorphism of $\{0,1,\ldots,n-1\}$-crystals.}
\end{split}
\end{equation*}
In this section we show in Theorem~\ref{thm:uniqueness} that, given $\Psi_0$ and $\Psi_n$,
there exists a unique $\{0,1,\ldots,n\}$-crystals isomorphism $\Psi:V^{r,s} \simeq B^{r,s}$
(with the exception of $r=n$ for $C_n^{(1)}$ and $D_{n+1}^{(2)}$ which is treated in 
Section~\ref{subsec:exceptional CD}).

We first prepare three preliminary lemmas that are used in the proof.

\begin{lemma}\label{lem:eps}
Let $b\in V^{r,s}$ be a $J$-lowest weight element. If $b$ contains a $\ol{p}$, with $p \neq 1$, then there 
exists a sequence $\vec{a}$ with letters in $\{2,\ldots,n-1 \}$, such that $\ve_n(e_{\vec{a}} b) > 0$.
\end{lemma}

\begin{proof}
By Lemma~\ref{lemma:Jlowest} the form of $b$ is known.
Let $\ell$ be the rightmost column, that contains a barred letter $\ol{p}$ which is not $\ol{1}$.
If there is no unbarred letter in column $\ell$, then by operating with
$e_{n-1}^{\max} \cdots e_p^{\max}$ one changes the $\ol{p}$ to $\ol{n}$. This
cannot be bracketed below, since there is no unbarred letter below and to the right, so
$\ve_n(e_{n-1}^{\max} \cdots e_p^{\max} b) > 0$.
Otherwise column $\ell$ has the form
\begin{equation*}
\begin{array}{|c|}
\hline
\ol{l}_t^{\phantom{L}}\\ \hline
\vdots\\ \hline
\ol{l}_1^{\phantom{L}}\\ \hline
n\\ \hline
n-1\\ \hline
\vdots\\ \hline
k \\ \hline
\end{array}
\qquad \text{or} \qquad
\begin{array}{|c|}
\hline
\ol{l}_t^{\phantom{L}}\\ \hline
\vdots\\ \hline
\ol{l}_1^{\phantom{L}}\\ \hline
0^\alpha\\ \hline
n\\ \hline
\vdots\\ \hline
k \\ \hline
\end{array}
\end{equation*}
where $l_i<k$ for all $1\le i\le t$, $0^\alpha$ stands for $\alpha$ boxes filled with $0$, and to the right 
of this column there are only unbarred letters and possibly $\ol{1}$. By operating with 
$e_{k-2}^{\max} \cdots e_{l_1}^{\max}$ the rightmost box that is changed is the $\ol{l}_1$ of column
$\ell$. Applying in addition $e_{n-2}^{\max} \cdots e_{k-1}^{\max}$ column $\ell$ becomes
\begin{equation*}
\begin{array}{|c|}
\hline
\ol{l}_t^{\phantom{L}}\\ \hline
\vdots\\ \hline
\ol{n-1}^{\phantom{L}}\\ \hline
n\\ \hline
n-2\\ \hline
\vdots\\ \hline
k-1 \\ \hline
\end{array}
\qquad \text{or} \qquad
\begin{array}{|c|}
\hline
\ol{l}_t^{\phantom{L}}\\ \hline
\vdots\\ \hline
\ol{n-1}^{\phantom{L}}\\ \hline
0^\alpha\\ \hline
n\\ \hline
n-2\\ \hline
\vdots\\ \hline
k-1 \\ \hline
\end{array}
\end{equation*}
and the columns to the right of column $\ell$ become
\begin{equation*}
\begin{array}{|c|}
\hline
\ol{1}^{\phantom{L}}\\ \hline
0^\beta\\ \hline
n\\ \hline
n-2\\ \hline
\vdots\\ \hline
m \\ \hline
\end{array}
\qquad \text{or} \qquad
\begin{array}{|c|}
\hline
0^\beta\\ \hline
n\\ \hline
n-2\\ \hline
\vdots\\ \hline
m \\ \hline
\end{array}\, .
\end{equation*}
where $m \geq k-1$ and $\beta\ge 0$.
By operating with $e_{n-1}^{\max}$, the letter $\ol{n-1}$ in column $\ell$ becomes $\ol{n}$ and all 
the $n$ below and to the right become $n-1$. If there is a $0$ it is unbracketed below. If there is no
$0$, then $\ol{n}$ is unbracketed below. Therefore $\ve_n(e_{n-1}^{\max} \cdots e_{l_1}^{\max} b) > 0$.
\end{proof}

\begin{lemma} 
\label{lem:not_height_n-1}
Let $b \in V^{n,s}$ be $J$-lowest such that $b$ does not contain any $\ol{k}$ with $k \neq 1$,
contains at least one $\ol{1}$ at height $n$, and has at least one column of height strictly
less than $n-1$.
Then there exists a sequence $\vec{a}$ with elements in $\{1,\ldots, n-1\}$
such that $\ve_n (e_{\vec{a}} b) > 0$ and $\vp_0 (e_{\vec{a}} b) > 0$.
\end{lemma}

\begin{proof}
Since $b$ is $J$-lowest, its columns of height $n$ from top to bottom must be
$n,n-1,\ldots,1$ or $\ol{1}, n, n-1, \ldots, 2$ or for the classical subalgebra $B_n$ possibly
$0^\alpha, n, n-1,\ldots, 1+\alpha$ with $\alpha\ge 1$. Since no more than one $0$ can appear in 
the same row, there can be at most one column of the form $0^\alpha, n, n-1,\ldots, 1+\alpha$.
Furthermore, since $b$ is supposed to contain at least one $\ol{1}$ at height $n$,
we must have $\alpha=1$ if a column containing $0$ exists.

Set $k=2\ell+j$ if there is no column of height $n$ containing $0$ and 
$k=2\ell+j-1$ otherwise, where $\ell$ is the number of columns of height $n$ not of the form
$n,n-1,\ldots,1$ and $j$ is the number of columns of $b$ of height $n-1$ of the form
$n,n-1,\ldots, 2$. Then $e_1^k(b)$ is a tableau, where the columns of height
$n$ not of the form $n,n-1,\ldots,1$ consist of the letters $\ol{2}, n, n-1, \ldots, 3, 1$ or 
$0, n, n-1, \ldots, 3, 1$. By passing to the $\{2,\ldots, n-1\}$-highest weight element, one obtains a 
tableau $e_{\vec{a}}b$ without $n$ (except in the columns of height $n$ of the form
$n,n-1,\ldots,1$), but with $\ol{n}$ or $0$ in the last $\ell$ columns of height $n$, 
where $\vec{a}$ has elements in $\{1,\ldots, n-1\}$.
Hence $\ve_n(e_{\vec{a}}b) > 0$. The corresponding $\pm$-diagram possibly has $+$ signs in 
rows of height $n$ or $n-1$, but no other $+$ signs. Since by assumption
there is at least one column of height strictly less than $n-1$, it follows 
from Proposition~\ref{prop:phi0} (2) that $\varphi_0(e_{\vec{a}}b) > 0$.
\end{proof}

\begin{lemma} \label{lem:not_related}
Let $b_1,b_2\in V^{r,s}$ be $J$-lowest, $b_1\neq b_2$, $\wt(b_1)=\wt(b_2)$, $b_1 \sim_{J'} b_2$,
and assume that $b_1$ does not contain any $\ol{k}\neq \ol{1}$.
If $b_1, b_2$ differ in boxes in rows strictly below $n-1$, we have
\begin{equation*}
\begin{split}
 	&\varphi_0(e_1^{\max}b_1), \varphi_0(e_1^{\max}b_2)>0 \quad \text{and}\\
	&e_1^{\max} b_1 \not \sim_{J'} e_1^{\max} b_2.
\end{split}
\end{equation*}
\end{lemma}

\begin{proof}
By Corollary~\ref{cor:J-lt}, $b_1$ and $b_2$ have the same inner tableau and only differ in 
the positions of the $\ol{1}$; the number of $\ol{1}$ in $b_1$ and $b_2$ must be the same 
since $\wt(b_1)=\wt(b_2)$.

By Lemma~\ref{lemma:Jlowest}, $b_1$ and $b_2$ contain $\ol{1}$ in certain positions
and the remainder of each column is filled with $n,n-1,\ldots,k$ for some $k$.
The $\pm$-diagrams corresponding to $b_1$ and $b_2$ contain columns with $-$ 
(in the position where $b_1$ and $b_2$ contain $\ol{1}$, respectively) and
otherwise only empty columns. By acting with $e_1^{\max}$ every $2$ changes into
$1$ and every $\ol{1}$ changes to $\ol{2}$. The $\pm$-diagram of $e_1^{\max}b_i$ will 
possibly have some $+$ at height $n$ and $n-1$, and otherwise only contains empty columns.
Therefore the inner shape of $e_1^{\max}b_1$ and $e_1^{\max}b_2$ is different
and hence $e_1^{\max} b_1 \not \sim_{J'} e_1^{\max} b_2$. By Proposition~\ref{prop:phi0} (2)
we conclude that  $\varphi_0(e_1^{\max}b_1), \varphi_0(e_1^{\max}b_2)>0$.
\end{proof}

\begin{theorem} 
\label{thm:uniqueness}
Let $1\le r\le n$ for type $A_{2n}^{(2)}$, and $1\le r<n$ for types $C_n^{(1)}$ and $D_{n+1}^{(2)}$.
We have $\Psi_0(b)=\Psi_n(b)$ for all $b\in V^{r,s}$ and hence there exists a unique $I$-crystal 
isomorphism $\Psi:V^{r,s}\cong B^{r,s}$.
\end{theorem}

\begin{proof}
Since $e_i \Psi_0=\Psi_0 e_i$ and $e_i \Psi_n = \Psi_n e_i$ for $i\in J$,
it suffices to show that $\Psi_n(b)=\Psi_0(b)$ for $J$-lowest elements $b$. We prove this
by downward induction on $(\wt b,\sum_{i=1}^n\epsilon_i)$, where $\epsilon_i$ are the
canonical basis vectors in $P=\Z^n$ and $(\cdot,\cdot)$ is the canonical inner product;
the quantity $(\wt b,\sum_{i=1}^n\epsilon_i)$ corresponds to the difference between the number
of unbarred and barred letters in the tableau of $b$.
If this value is maximal, there is only one $J$-lowest element.

First assume that $b$ is a $J$-lowest element satisfying
\begin{quote}
  $\vp_0(b)>0$ and
 \begin{enumerate}
	\item $b$ contains $\ol{k}$ with $k>1$ or 
	\item $b$ satisfies the conditions of Lemma~\ref{lem:not_height_n-1}.
\end{enumerate}
\end{quote}
In case 1 by Lemma \ref{lem:eps} there exists a sequence $\vec{a}$ consisting of elements in
$\{2,\ldots,n-1\}$ such that $\ve_n(e_{\vec{a}}b)>0$. Then we also have $\vp_0(e_{\vec{a}}b)>0$.
In case 2 by Lemma~\ref{lem:not_height_n-1} there exists a sequence $\vec{a}$ with elements
in $\{1,2,\ldots, n-1\}$ such that $\vp_0(e_{\vec{a}}b)>0$ and $\ve_n(e_{\vec{a}}b)>0$. Set 
$b'=e_{\vec{a}}b$. We have
\begin{align*}
e_nf_0\Psi_n(b')&=e_n\Psi_n(f_0b')=e_n\Psi_0(f_0b')=\Psi_0(e_nf_0b')\\
&=\Psi_n(f_0e_nb')=f_0\Psi_n(e_nb')=f_0\Psi_0(e_nb')=f_0e_n\Psi_0(b').
\end{align*}
Here in the 2nd, 4th and 6th equality we have used the induction hypothesis. Hence we have
$\Psi_n(b)=\Psi_0(b)$ in this case.

Next assume that the $J$-lowest element $b$ satisfies
\begin{quote}
  $\vp_0(b)>0$ and \newline
  $b$ does not contain $\ol{k}$ with $k>1$ and \newline
  $b$ does not satisfy the conditions of Lemma~\ref{lem:not_height_n-1}.
\end{quote}
Suppose
\begin{equation} \label{eq2}
\Psi_n(b)=\Psi_0(b')
\end{equation}
for $b,b'$ such that $b\neq b'$. One can assume $b'$ is also $J$-lowest and $\wt b=\wt b'$.
We show by contradiction that this is not possible. We have
\[
\Psi_0(b')=e_0\Psi_n(f_0b)=e_0\Psi_0(f_0b).
\]
The second equality is due to the induction hypothesis. From the equality of the LHS and RHS
we have $b\sim_{J'}b'$. If $b'$ contains $\ol{k}$ with $k>1$ or $b'$ satisfies the conditions of
Lemma~\ref{lem:not_height_n-1}, we already know $\Psi_n(b')=\Psi_0(b')$
from the previous case. Hence we can assume that $b'$ does not contain $\ol{k}$ with $k>1$ 
and does not satisfy the conditions of Lemma~\ref{lem:not_height_n-1} either.

Suppose that $b$ and $b'$ do not differ in boxes in rows strictly below $n-1$. Since $b\sim_{J'} b'$,
this means by Corollary~\ref{cor:J-lt} that they differ in boxes containing $\ol{1}$ in row $n$ and
$n-1$. But then at least one of $b$ or $b'$ satisfies the conditions of Lemma~\ref{lem:not_height_n-1}
which is a contradiction. Hence $b$ and $b'$ must differ in boxes in row strictly below $n-1$.
Then by Lemma~\ref{lem:not_related} we have $\vp_0(e_1^{\max}b)>0$ and
\begin{equation} \label{not related}
e_1^{\max}b\not\sim_{J'}e_1^{\max}b'.
\end{equation}
{}From \eqref{eq2} one has
\[
\Psi_0(e_1^{\max}b')=\Psi_n(e_1^{\max}b)=e_0\Psi_n(f_0e_1^{\max}b)=e_0\Psi_0(f_0e_1^{\max}b).
\]
In the last equality we used the induction hypothesis. But the equality of the LHS and RHS
contradicts \eqref{not related}.

We are left to show $\Psi_n(b)=\Psi_0(b)$ when $\vp_0(b)=0$. However, from 
Lemma~\ref{lem:decomp1} for type $C_n^{(1)}$ and Lemma~\ref{lem:decomp2} for types 
$A_{2n}^{(2)}$ and $D_{n+1}^{(2)}$ such an element $b$ is unique if we specify the weight. So one 
has to have $\Psi_n(b)=\Psi_0(b)$ also in this case. This completes the proof.
\end{proof}

\section{KR crystals for exceptional nodes}
\label{sec:exceptional}

A node $r$ in the Dynkin diagram is called special if there is a Dynkin diagram automorphism
that maps $r$ to $0$. When $r$ is a special node, the corresponding KR crystal $B^{r,s}$
is irreducible as a $\{1,2,\ldots,n\}$-crystal. For type $A_n^{(1)}$, all nodes $r$ are special and have 
already been treated in Section~\ref{subsec:A}. For types $B_n^{(1)}$, $A_{2n-1}^{(2)}$ and
$D_n^{(1)}$, the node $r=1$ is special and has already been treated in Section~\ref{subsec:DBA}.
The remaining special nodes are $r=n$ for types $C_n^{(1)}$ and $D_{n+1}^{(2)}$, and 
$r=n-1,n$ for type $D_n^{(1)}$, for which the usual decomposition~\eqref{eq:decomp W} does
not hold anymore. The KR crystals for these exceptional nodes are treated in Sections~\ref{subsec:exceptional CD}
and~\ref{subsec:exceptional D}, respectively.

\subsection{$B^{n,s}$ of type $C_n^{(1)}, D_{n+1}^{(2)}$}
\label{subsec:exceptional CD}
In this section we give the combinatorial description of the KR-crystal $B^{n,s}$
of types $C_n^{(1)}$ and $D_{n+1}^{(2)}$. As a $\{1,2,\ldots,n\}$-crystal we have
the isomorphism
\begin{equation} \label{eq:C spin decomp}
	B^{n,s} \cong B(s\La_n).
\end{equation}

First consider type $C_n^{(1)}$. The elements in $B(s\La_n)$ are KN-tableaux of shape 
$(s^n)$ (see Section~\ref{subsec:KN C}). Recall from Section~\ref{subsec:branching}, 
that the $J'$-highest weight elements of shape $(s^n)$ are in bijection with $\pm$-diagrams.
Since all columns are of height $n$ and the classical subalgebra of $C_n^{(1)}$ is $C_n$, each 
column is either filled with $+$, $-$, or $\mp$.
Hence, if there are $\ell_1$ columns containing $+$, $\ell_2$ columns containing $-$, and 
$\ell_3$ columns containing $\mp$, we may identify $\pm$-diagrams $P$ with triples 
$(\ell_1, \ell_2, \ell_3)$ such that $\ell_1 + \ell_2 + \ell_3 = s$ and $\ell_1,\ell_2,\ell_3\ge 0$. 

In order to describe the affine structure, it suffices to define $e_0$ on such triples, since 
$e_0$ commutes with $e_2,\ldots, e_n$.
Acting with $e_0$ changes neither the inner shape of $P$ (since $e_0$ commutes with 
$e_2,\ldots, e_n$) nor the outer shape of $P$ (by the decomposition~\eqref{eq:C spin decomp}).
Hence $\ell_3$ is invariant under $e_0$. 
If $\ve_0(\ell_1, \ell_2, \ell_3) > 0$, for weight reasons we must have $e_0(\ell_1, \ell_2, \ell_3) 
= (\ell_1 - 1, \ell_2 + 1, \ell_3)$. Again for weight reasons, 
if $\ell_1 = 0$, then $e_0(\ell_1, \ell_2, \ell_3) = \emptyset$. We now calculate 
$\ve_0 (\ell_1, \ell_2, \ell_3)$.

\begin{lemma}
$\ve_0(\ell_1,\ell_2,\ell_3) = \ell_1$.
\end{lemma}
\begin{proof}
If $\ell_2 = 0$, then $\ve_0(\ell_1,0,\ell_3) -\vp_0(\ell_1, 0, \ell_3) = \ell_1$ and therefore 
$\ve_0(\ell_1, 0, \ell_3) \geq \ell_1$. But by the previous observation 
$\ve_0(\ell_1, \ell_2, \ell_3) \leq \ell_1$. The claim follows.
\end{proof}

\begin{definition}
\label{def:C exceptional}
The combinatorial crystal $V^{n,s}$ of type $C_n^{(1)}$ is defined to be $B(s\La_n)$ as 
a $\{1,2,\ldots,n\}$-crystal. The action of $f_0$ and  $e_0$ on $\{2,3,\ldots,n\}$-highest weight 
elements is given by
\begin{equation*}
\begin{split}
	f_0(\ell_1,\ell_2,\ell_3) &= \begin{cases}
	(\ell_1 + 1, \ell_2 - 1, \ell_3) & \text{if $\ell_2>0$,}\\
	\emptyset & \text{otherwise,}
	\end{cases}\\
	e_0(\ell_1,\ell_2,\ell_3) &= \begin{cases}
	(\ell_1 - 1, \ell_2 + 1, \ell_3) & \text{if $\ell_1>0$,}\\
	\emptyset & \text{otherwise.}
	\end{cases}
\end{split}
\end{equation*}
\end{definition}

Next consider type $D_{n+1}^{(2)}$ whose classical subalgebra is of type $B_n$.
Since $\La_n=\frac{1}{2}(\epsilon_1+\cdots+\epsilon_n)$, the elements in $B(s\La_n)$
are KN-tableaux of shape $((s/2)^n)$ when $s$ is even and of shape $(((s-1)/2)^n)$ plus an extra
spin column when $s$ is odd. By Section~\ref{subsec:branching}, the $J'$-highest weight
elements are in bijection with $\pm$-diagrams, where columns of height $n$ can contain
$+$, $-$, $\mp$ and at most one $0$; the spin column of half width can contain $+$ or $-$.

We may again encode a $\pm$-diagram $P$ as a triple $(\ell_1, \ell_2, \ell_3)$,
where $\ell_1$ is twice the number of columns containing a single $+$ sign, $\ell_2$ is twice 
the number of columns containing a single $-$ sign (where spin column are counted as $1/2$ 
columns), and $\ell_3$ is twice the number of columns containing $\mp$.
If $P$ contains a $0$-column, then $\ell_1 + \ell_2 + \ell_3 =s-2$, otherwise  
$\ell_1 + \ell_2 + \ell_3 =s$.

As in the case $C_n^{(1)}$, since $e_0$ commutes with $e_2,\ldots,e_n$ it suffices to specify
the action of $e_0$ on $\{2,3,\ldots,n\}$-components or equivalently on triples $(\ell_1,\ell_2,\ell_3)$.

\begin{lemma}
$\ve_0(\ell_1, \ell_2, \ell_2) = \ell_1 + \gamma$, where $\gamma$ is $1$ if there is a $0$ column, 
and $0$ otherwise.
\end{lemma}
\begin{proof}
The proof is the same as above for $C_n^{(1)}$.
\end{proof}

\begin{definition}
\label{def:D2 exceptional}
The combinatorial crystal $V^{n,s}$ of type $D_{n+1}^{(2)}$ is defined to be $B(s\La_n)$ as 
a $\{1,2,\ldots,n\}$-crystal. The action of $f_0$ and  $e_0$ on $\{2,3,\ldots,n\}$-highest weight 
elements is given by
\begin{equation*}
\begin{split}
	f_0(\ell_1,\ell_2,\ell_3) &= \begin{cases}
	(\ell_1+2, \ell_2, \ell_3) & \text{if $\ell_1 + \ell_2 + \ell_3 < s$,}\\
	(\ell_1, \ell_2-2, \ell_3) & \text{if $\ell_1 + \ell_2 + \ell_3 = s$ and $\ell_2 > 1$,}\\
	(\ell_1+1, 0, \ell_3) & \text{if $\ell_1 + \ell_2 + \ell_3 = s$ and $\ell_2 = 1$,}\\
	\emptyset & \text{if $\ell_1 + \ell_2 + \ell_3 = s$ and $\ell_2=0$,}
	\end{cases}\\
	e_0(\ell_1,\ell_2,\ell_3) &= \begin{cases}
	(\ell_1, \ell_2 + 2, \ell_3) & \text{if $\ell_1 + \ell_2 + \ell_3 < s$,}\\
	(\ell_1-2, \ell_2, \ell_3) & \text{if $\ell_1 + \ell_2 + \ell_3 = s$ and $\ell_1 > 1$,}\\
	(0, \ell_2+1, \ell_3) & \text{if $\ell_1 + \ell_2 + \ell_3 = s$ and $\ell_1 = 1$,}\\
	\emptyset & \text{if $\ell_1 + \ell_2 + \ell_3 = s$ and $\ell_1=0$.}
	\end{cases}
\end{split}
\end{equation*}
\end{definition}

\begin{theorem} \label{thm:exceptional CD}
For types $C_n^{(1)}$ and $D_{n+1}^{(2)}$, we have
\begin{equation*}
	V^{n,s} \cong B^{n,s}.
\end{equation*}
\end{theorem}

\begin{proof}
By~\cite[Theorem 1.1]{OS:2008} the KR crystal $B^{n,s}$ exists. Given the 
decomposition~\eqref{eq:C spin decomp} as a $\{1,2,\ldots,n\}$-crystal, the construction of 
$V^{n,s}$ as an affine crystal was uniquely specified by weight considerations. Hence, since
$V^{n,s}$ and $B^{n,s}$ have the same $\{1,2,\ldots,n\}$-decomposition and $B^{n,s}$ exists,
$V^{n,s}$ and $B^{n,s}$ must be isomorphic.
\end{proof}

\subsection{$B^{n,s}$ and $B^{n-1,s}$ of type $D_n^{(1)}$}
\label{subsec:exceptional D}
In this section we give a combinatorial model for the KR-crystals $B^{n,s}$ and $B^{n-1,s}$ of 
type $D_n^{(1)}$ that are associated to the spin nodes $n-1$ and $n$ in the Dynkin diagram.
As $\{1,2,\ldots, n\}$-crystals we have the isomorphisms
\begin{equation} \label{eq:spin D decomp}
\begin{split}
	B^{n,s} &\cong B(s\La_n),\\
	B^{n-1,s} & \cong B(s\La_{n-1}).
\end{split}
\end{equation}
The combinatorial KR-crystals $V^{n,s}$ and $V^{n-1,s}$ are constructed to have the same 
classical decomposition as in~\eqref{eq:spin D decomp}. To define the affine crystal action,
we first introduce an involution $\sigma:B^{n,s} \leftrightarrow B^{n-1,s}$ corresponding to the 
Dynkin diagram automorphism that interchanges the nodes $0$ and $1$.
Under this involution, $\{2,3,\ldots,n\}$-components need to be mapped to 
$\{2,3,\ldots,n\}$-components. Hence it suffices to define $\sigma$ on 
$\{2,3,\ldots,n\}$-highest weight elements or equivalently $\pm$-diagrams.
Recall from Section~\ref{subsec:branching}, that for weights $\La=s\La_n$ or $s\La_{n-1}$,
the $\pm$-diagram can contain columns with $+$ and $\mp$ or with $-$ and $\mp$ (but not a 
mix of $-$ and $+$ columns).

\begin{definition}
\label{def:sigma exceptional}
The involution $\sigma:B^{n,s} \leftrightarrow B^{n-1,s}$ maps a $\pm$-diagram $P$
to a $\pm$-diagram $P'$ of opposite color where columns containing $+$ are interchanged
with columns containing $-$ and vice versa.
\end{definition}

\begin{definition}
\label{def:D exceptional}
The combinatorial crystal $V^{n,s}$ (resp. $V^{n-1,s}$) of type $D_n^{(1)}$ is defined to be 
$B(s\La_n)$ (resp. $B(s\La_{n-1})$) as a $\{1,2,\ldots,n\}$-crystal. The action of $f_0$ and  $e_0$ is
\begin{equation*}
	e_0 = \sigma \circ e_1 \circ \sigma \qquad \text{and} \qquad f_0 = \sigma \circ f_1 \circ \sigma
\end{equation*}
with $\sigma$ as in Definition~\ref{def:sigma exceptional}.
\end{definition}

\begin{theorem} \label{thm:exceptional D}
For type $D_n^{(1)}$, we have 
\begin{equation*}
\begin{split}
	V^{n,s} &\cong B^{n,s},\\
	V^{n-1,s} &\cong B^{n-1,s}.
\end{split}
\end{equation*}
\end{theorem}

\begin{proof}
$V^{n,s}$ and $B^{n,s}$ have the same decomposition as $\{1,2,\ldots,n\}$ and 
$\{0,2,\ldots,n\}$-crystals
\begin{align*}
	\Psi_0: & \quad V^{n,s} \cong B^{n,s} \cong B(s\La_n) 
	&& \text{as a $\{1,2,\ldots,n\}$-crystals,}\\
	\Psi_1: & \quad V^{n,s} \cong B^{n,s} \cong B(s\La_{n-1}) 
	&& \text{as a $\{0,2,\ldots,n\}$-crystals.}
\end{align*}
The $\{1,2,\ldots,n\}$-crystal isomorphism $V^{n,s}\cong B(s\La_n)$ is true
by definition and the $\{0,2,\ldots,n\}$-crystal isomorphism $V^{n,s}\cong B(s\La_{n-1})$ 
follows by the application of $\sigma$ since $e_0 = \sigma e_1 \sigma$ and 
$e_i = \sigma e_i \sigma$ for $i\neq 0,1$. Chari~\cite{Chari:2001} proved that
$B^{n,s}\cong B(s\La_n)$ as $\{1,2,\ldots,n\}$-crystals.
For the proof that $B^{n,s} \cong B(s\La_{n-1})$ as a
$\{0,2,\ldots,n\}$-crystal, it suffices to show there exists a
corresponding highest weight vector, since the crystal is irreducible.
Applying the Weyl group element $r_\beta$, which is the reflection for the root
$\beta=\epsilon_1+\epsilon_n$, to the $\{1,2,\ldots,n\}$-highest weight element
yields the $\{0,2,\ldots,n\}$-highest weight element.

Note that if $\Psi_0(b)=\Psi_1(b)$ for a $b$ in a given $D_{n-1}$-component $\mathcal{C}$, then 
$\Psi_0(b')=\Psi_1(b')$ for all $b'\in \mathcal{C}$ since $e_i\Psi_0(b')=\Psi_0(e_ib')$ and $e_i\Psi_1(b')=
\Psi_1(e_ib')$ for $i\in J'=\{2,3,\ldots,n\}$.
Furthermore observe that by Remark~\ref{rem:psi weight} $\Psi_0$ and $\Psi_1$ preserve weights, that is,
$\wt(b)=\wt(\Psi_0(b))=\wt(\Psi_1(b))$ for all $b\in V^{n,s}$.

Since $e_i$ commutes with $\Psi_0$ and $\Psi_1$ for $i\in J'$, it follows that $J'$-components
in $V^{n,s}$ must map to $J'$-components in $B^{n,s}$. However, as can be seen from the 
description using $\pm$-diagrams, the branching $D_n\to D_{n-1}$ on $B(s\La_n)$ and 
$B(s\La_{n-1})$ is multiplicity-free once the weight is fixed. Hence we must have 
$\Psi_0(b)=\Psi_1(b)$ for all $b\in V^{n,s}$. The proof for $V^{n-1,s}$ is analogous.
\end{proof}

\section{Dynkin automorphism for type $C_n^{(1)}$ and $D_{n+1}^{(2)}$}
\label{sec:auto}

By construction, the Dynkin diagram automorphism for type $A_{n-1}^{(1)}$, $B_n^{(1)}$, 
$D_n^{(1)}$, and $A_{2n-1}^{(2)}$ acts on the combinatorial crystal $V^{r,s}$, except for
$r=n-1,n$ for type $D_n^{(1)}$.
The Dynkin diagrams for type $C_n^{(1)}$ and $D_{n+1}^{(2)}$ also have an automorphism
mapping $i\mapsto n-i$ for all $i\in \{0,1,\ldots,n\}$. However, from the construction of $V^{r,s}$
for these types using Dynkin diagram foldings and similarity methods, it is not obvious that
this Dynkin diagram automorphism extends to $V^{r,s}$. This is proven in Theorem~\ref{thm:auto}.
This shows in particular that~\cite[Assumption 1]{FSS:2007} holds, which was used to show 
that the classical isomorphism from the Demazure crystal to the KR crystal,
sends zero arrows to zero arrows. 

\begin{theorem} \label{thm:auto}
Let $B^{r,s}$ ($1\le r\le n,s\ge1$) be the KR crystal for type $C_n^{(1)},D_{n+1}^{(2)}$.
Then there exists an involution $\sigma$ on $B^{r,s}$ satisfying 
\begin{equation} \label{i to n-i}
\sigma\circ e_i=e_{n-i}\circ\sigma \quad \text{and} \quad \sigma\circ f_i=f_{n-i}\circ\sigma
\quad \text{for all  $i\in I$.}
\end{equation}
\end{theorem}

Let $\widetilde{\mathcal{B}}^{r,s}$ be the fixed point crystal in Case (a) of Section 2.1
constructed in~\cite{NS:2006}. Naito and Sagaki~\cite{NS:2006} show that 
$\widetilde{\mathcal{B}}^{r,s}$ is a regular $D_{n+1}^{(2)}$-crystal which decomposes into 
$\bigoplus_\La B(\La)$ for $1\le r<n$ and $B(s\La_n)$ for $r=n$ as a $B_n$-crystal, where the 
direct sum $\bigoplus_\La$ is over all $\La$ obtained from $s\La_r$ by removing boxes and $B(\La)$ 
is a highest weight $B_n$-crystal of highest weight $\La$. They also show that as a 
$D_{n+1}^{(2)}$-crystal it is isomorphic to the virtual $U'_q(D_{n+1}^{(2)})$-crystal defined 
in~\cite[Section 6.7]{OSS:2003}.

The next lemma is used for the proof of type $D_{n+1}^{(2)}$.

\begin{lemma} \label{lem:auto D}\mbox{}
\begin{itemize}
\item[(1)] The decomposition
\[
\widetilde{\mathcal{B}}^{r,s}\simeq 
\begin{cases}
\bigoplus_\La B(\La)&\text{ for }1\le r<n\\
B(s\La_n)&\text{ for $r=n$}
\end{cases}
\]
also holds as a $\{0,1,\ldots,n-1\}$-crystal. 
\item[(2)] There exists an involution $\sigma$ on $\widetilde{\mathcal{B}}^{r,s}$ satisfying \eqref{i to n-i}.
\end{itemize}
\end{lemma}

\begin{proof}
Note that the Weyl group of type $B_n$ contains an element that maps 
$\La_j-2\La_0$ to $\La_{n-j}-2\La_n$ for all $1\le j\le n$. Hence, the decomposition as a 
$\{0,1,\ldots,n-1\}$-crystal follows from that as a $\{1,2,\ldots,n\}$-crystal as in the proof of 
Lemma~\ref{lem:decomp1}.

We prove (2). Set $\sigma=\pr^n$, where $\pr$ is the promotion operator on the ambient 
$A_{2n-1}^{(1)}$-crystal as defined in Section~\ref{subsec:A}. The map $\sigma$ 
satisfies~\eqref{i to n-i} with $x_i=\hat{x}_i$ for $i=0,n$ and $x_i=\hat{x}_i\hat{x}_{2n-i}$ otherwise,
where $x=e,f$. Hence, it suffices to show
\begin{equation} \label{sigma incl}
\sigma(\widetilde{\mathcal{B}}^{r,s})\subset\widetilde{\mathcal{B}}^{r,s}.
\end{equation}
Let $u$ be the unique dominant extremal element in the ambient $A_{2n-1}^{(1)}$-crystal.
Then $u$ belongs to $\widetilde{\mathcal{B}}^{r,s}$. The inclusion \eqref{sigma incl} is now clear, since 
$\widetilde{\mathcal{B}}^{r,s}$ is generated from $u$ by applying $e_i$ and $f_i$.
\end{proof}

\begin{proof}[Proof of Theorem~\ref{thm:auto} for type $D_{n+1}^{(2)}$]
By Theorem \ref{thm:uniqueness} we know $B^{r,s}\simeq \widetilde{\mathcal{B}}^{r,s}$ as 
$D_{n+1}^{(2)}$-crystals. Hence $\sigma$ exists by (2) of  Lemma~\ref{lem:auto D}.
\end{proof}

We prepare a lemma and a proposition for the proof of type $C_n^{(1)}$.

\begin{lemma} \label{lem:same prop}
Let $V^{r,s}_{C_n^{(1)}}$ be the $C_n^{(1)}$-crystal constructed in Section \ref{subsec:C}.
Let $b$ be a $\{2,\ldots,n\}$-highest element of $V^{r,2s}_{C_n^{(1)}}$ whose shape is obtained from 
an $r\times2s$ rectangle by removing $1\times4$ rectangular pieces. Let $P$ be the corresponding 
$\pm$-diagram of $b$. Suppose that the number of columns with $\mp,+,-,\cdot$ at each height is
even. Then the shape of $e_0^2(b)$ and the corresponding $\pm$-diagram have the same property.
\end{lemma}

\begin{proof}
Let $\Vh^{r,s}$ be the ambient crystal of type $A_{2n-1}^{(2)}$ in the definition of the type $C_n^{(1)}$
KR crystal. We follow the same set-up as in the proof of Lemma~\ref{lem:double}. For elements
in $b\in V^{r,s}:=V^{r,s}_{C_n^{(1)}}$, we have $\sigma(b)=b$. Hence on $V^{r,s}$ we have
$e_0^2 = \eh_1 \eh_0 \eh_1 \eh_0= \eh_1^2 \eh_0^2 = \eh_1^2 (\sigma \eh_1 \sigma) (\sigma \eh_1
\sigma) = \eh_1^2 \sigma \eh_1^2$. The proof follows that of Lemma~\ref{lem:double} with
all columns and operations doubled.
\end{proof}

\begin{proposition} \label{prop:auto C}\mbox{}
\begin{itemize}
\item[(1)] There exists a $C_n^{(1)}$-crystal $\widetilde{V}^{r,s}$ and a unique injective map 
	$S:\widetilde{V}^{r,s}\longrightarrow V^{r,s}_{D_{n+1}^{(2)}}$ such that 
\[
S(e_ib)=e_i^{m_i}S(b),\quad S(f_ib)=f_i^{m_i}S(b) \text{ for }i\in I,
\]
	where $(m_0,m_1,\ldots,m_{n-1},m_n)=(2,1,\ldots,1,2)$.

\item[(2)] $\widetilde{V}^{r,s}$ is connected.

\item[(3)] There exists an involution $\sigma$ on $\widetilde{V}^{r,s}$ satisfying \eqref{i to n-i}.
\end{itemize}
\end{proposition}

\begin{proof}
We only prove the case $1\le r\le n-1$. The case $r=n$ is similar and easier.

Let us prove (1). Let $B(\La)$ (resp. $B_{B_n}(\La)$) 
be the $C_n$ (resp. $B_n$)-crystal of the highest weight module of highest weight $\La$. 
By~\cite[Theorem~5.1]{Ka:1996}, there exists a unique injective map 
$\ol{S}:B(\La)\rightarrow B_{B_n}(\La)$ such that $\ol{S}(e_ib)=e_i^{m_i}\ol{S}(b)$ and 
$\ol{S}(f_ib)=f_i^{m_i}\ol{S}(b)$ for $i\in I_0=\{1,2,\ldots,n\}$.

Now define $\widetilde{V}^{r,s}$, as a $C_n$-crystal, by
\begin{equation} \label{eq:Vtilde decomp}
   \widetilde{V}^{r,s}=\bigoplus_{\La}B(\La),
\end{equation}
where the sum is over all $\La$ obtained from $s\La_r$ by removing horizontal dominoes. From 
the above explanation we have an injective map
$S:\widetilde{V}^{r,s}\rightarrow V^{r,s}_{D_{n+1}^{(2)}}$. We introduce the 0-action on 
$\widetilde{V}^{r,s}$. 
Recall the construction of an injective map $S$ (denoted here by $S'$) $:V^{r,s}_{D_{n+1}^{(2)}}
\rightarrow V^{r,2s}_{C_n^{(1)}}$ of Section~\ref{subsec:A(2)D(2)}. For a $\{2,\ldots,n\}$-highest 
element $b$ in $V^{r,s}_{D_{n+1}^{(2)}}$, it is shown that $S'S(b)$ satisfies the assumptions of 
Lemma~\ref{lem:same prop} by using Lemma~\ref{lem:double pm}. Hence, $e_0^2S'S(b)$ also has 
the same property, so there exists a $b'\in \widetilde{V}^{r,s}$ such that $S'S(b')=e_0^2S'S(b)$. One 
can define $e_0(b)=b'$ on $\widetilde{V}^{r,s}$
and we have $S(e_0b)=e_0^2S(b)$. The case $f_0$ is similar.

Next we prove (2). Suppose that $\widetilde{V}^{r,s}$ is not connected.
Then there must exist $\La$, $\La'$ in the decomposition~\eqref{eq:Vtilde decomp}
that differ only in one horizontal domino, but lie in different components (otherwise all $\La$
would lie in the same component and hence $\widetilde{V}^{r,s}$ would be connected).
By Lemma~\ref{lem:double} there exists a $\{2,3,\ldots,n\}$-highest weight element $b\in V^{r,s}$
of type $C_n^{(1)}$ such that $b\in B(\La)$ and $b':=e_0(b)\in B(\La')$. Identifying the classical
decomposition of $V^{r,s}$ as in Lemma~\ref{lem:decomp} and $\widetilde{V}^{r,s}$ as
in~\eqref{eq:Vtilde decomp}, we may consider the corresponding elements 
$\tilde{b}\in B(\La) \subset \widetilde{V}^{r,s}$ and $\tilde{b}'\in B(\La') \subset \widetilde{V}^{r,s}$.
Then $S'S(\tilde{b})\in B(2\La)$ and $S'S(\tilde{b}') \in B(2\La')$, and using Lemma~\ref{lem:double pm}
for the corresponding $\pm$-diagrams we find that $S'S(\tilde{b}') = e_0^2 S'S(\tilde{b})$.
This implies by the definition of $e_0$ on $\widetilde{V}^{r,s}$ that $\tilde{b}'= e_0(\tilde{b})$ 
which contradicts the assumption that $B(\La)$ and $B(\La')$ lie in different components of 
$\widetilde{V}^{r,s}$. Hence $\widetilde{V}^{r,s}$ must be connected.

Finally we prove (3). By (1) one can consider the problem in the image of $S$. By 
Theorem~\ref{thm:auto} for type $D_{n+1}^{(2)}$ we know that there exists an involution $\sigma$ 
satisfying \eqref{i to n-i} on $V^{r,s}_{D_{n+1}^{(2)}}$. Since $\text{Im}\;S$ is generated by 
$e_i^{m_i}$ using (2), it is clear that $\text{Im}\;S$ is closed under $\sigma$.
\end{proof}

\begin{proof}[Proof of Theorem~\ref{thm:auto} for type $C_n^{(1)}$]
By construction, $\widetilde{V}^{r,s}$ of Proposition \ref{prop:auto C} and $B^{r,s}$ have the 
same decomposition as $\{1,\ldots,n\}$-crystals. From Proposition~\ref{prop:auto C} (3) it follows
that $\widetilde{V}^{r,s}$ and $B^{r,s}$ have the same decomposition also as $\{0,\ldots,n-1\}$-crystals. 
Then by Theorem~\ref{thm:uniqueness} we know that $B^{r,s}\simeq \widetilde{V}^{r,s}$ as 
$C_n^{(1)}$-crystals. Hence $\sigma$ on $B^{r,s}$ exists by (3) of Proposition~\ref{prop:auto C}.
\end{proof}

\end{document}